\documentclass[10pt]{article}
\usepackage{latexsym,amssymb,dsfont,textcomp}
\usepackage[centertags]{amsmath}
\usepackage{fullpage}

\usepackage{amstext}
\def\R{\mathbb R}
\def\C{\mathbb C}
\def\N{\mathbb N}
\def\Z{\mathbb Z}

\newcommand{\I}{\mathbf{I}}

\newcommand{\tz}{{\tilde z}}
\newcommand{\tx}{{\tilde x}}
\newcommand{\ty}{{\tilde y}}
\newcommand{\la}{\label}
\newcommand{\beq}{\begin{equation}}
\newcommand{\eeq}{\end{equation}}
\newcommand{\eps}{\varepsilon}
\newcommand{\qed}{\protect~\protect\hfill $\Box$}

\begin{document}

\newtheorem{theorem}{Theorem}[section]
\renewcommand{\thetheorem}{\arabic{section}.\arabic{theorem}}
\renewcommand{\labelenumi}{\rm (\roman{enumi})}
\newtheorem{definition}[theorem]{Definition}
\newtheorem{deflem}[theorem]{Definition and Lemma}
\newtheorem{lemma}[theorem]{Lemma}
\newtheorem{example}[theorem]{Example}
\newtheorem{remark}[theorem]{Remark}
\newtheorem{remarks}[theorem]{Remarks}
\newtheorem{cor}[theorem]{Corollary}
\newtheorem{pro}[theorem]{Proposition}

\renewcommand{\theequation}{\thesection.\arabic{equation}}

\title{Nekhoroshev type stability results for Hamiltonian systems
       with an additional transversal component}
\author{{\sc Markus Kunze$^{1}$ \& David M.A.~Stuart$^{2}$} \\[2ex]
        $^{1}$ Universit\"at Duisburg-Essen, Fakult\"at f\"ur Mathematik, \\
        D\,-\,45117 Essen, Germany \\[1ex]
        $^{2}$ Centre for Mathematical Sciences, Wilberforce Road, \\
        Cambridge, CB3 OWA, England \\[2ex]
        {\bf Key words:} Hamiltonian systems, Nekhoroshev stability, normal forms}
\date{}
\maketitle
\begin{abstract}\noindent 
We prove exponential stability theorems of Nekhoroshev type
for motion in the neighbourhood of an elliptic fixed point in Hamiltonian
systems having an additional transverse component of arbitrary dimension. 
The estimates that we obtain are independent of this dimension. 
\end{abstract}

%%%%%%%%%%%%%%%%%%%%%%%%%%%%%%%%%%%%%%%%%%%%%%%%%%%%%%%%%%%%%%%%%%%%%%%%%%%%%%%%%%%%%%%%

\setcounter{equation}{0}

\section{Introduction and statement of main results}

An integrable Hamiltonian, written in action angle variables
$(I,\phi)=(I_1,\dots I_n,\phi_1,\dots\phi_n)\in\R^n\times(\R/\Z)^n$,
takes the form $H_{{\rm int}}(I_1,\dots I_n)$, and the corresponding
equations of motion imply that the action variables $I_j$ are constant
while the angle variables $\phi_j$ evolve at the uniform rate 
$\frac{\partial H_{{\rm int}}}{\partial I_j}$. For a nonintegrable perturbation 
of such a system, described by a smooth Hamiltonian of the form
$$
H_{{\rm int}}(I_1, \dots I_n)+\eps H_{{\rm pert}}(I_1,\dots I_n,\phi_1,\dots\phi_n)\,,
$$
Nekhoroshev proved the following {\em exponential stability} estimate
in \cite{nek77}: let $H_{{\rm int}}$ satisfy a condition known 
as steepness, then there exist positive numbers $R_0,T_0,\eps_0,a,b$
such that for all small $\eps$
\beq\la{nee}
|I(t)-I(0)|\leq R_0\eps^b\quad\hbox{for}\quad |t|\leq T_0
e^{(\frac{\eps_{\small 0}}{\eps})^a}\,.
\eeq
This says that for small $\eps$ the action variables are almost, or effectively,
constant since they vary little over exponentially long time scales. 
In fact the main theorem in \cite[\S4.4]{nek77} proves exponential stability
bounds for slightly more general perturbations 
$$H_{{\rm pert}}=H_{{\rm pert}}(I_1,\dots I_n,\phi_1,\dots\phi_n,\xi_1,
\dots\xi_N,\eta_1,\dots\eta_N)$$
in which there is dependence upon an additional set of $N$ (Darboux) conjugate
pairs $(\xi_j,\eta_j)$; we shall refer to these extra variables
as the {\em transverse component}.

It was also conjectured in \cite{nek77} that under appropriate
conditions such exponential stability should hold
in a sufficiently small neighbourhood of an elliptic equilibrium point.
Following a preliminary result in \cite[\S IV.2, Theorem 4]{lochak1}
this was proved in \cite{gfb,nied98} and then in \cite{poschel}
under convexity hypotheses which can be described as follows: let
\begin{itemize}
\item
$\{(x_j,y_j)\}_{j=1}^n$ be Darboux coordinates on $\R^{2n}$, and define
$I_j=(x_j^2+y_j^2)/2$,
\item
$\alpha\in\R^n$, and let $A$ be a {\em strictly positive }
$n\times n$ matrix, 
\item
$f$ be a real analytic function 
vanishing to fifth order at the origin,
\end{itemize}
then the dynamics in a neighbourhood of the origin in $\R^{2n}$
for the Hamiltonian
$$
H_0=\langle\alpha, I\rangle+\frac{1}{2}\,\langle AI, I\rangle+f
$$
satisfies exponential stability estimates; 
see theorem \ref{classicalthm} for a precise statement. 

In view of the above it is to be expected that exponential stability may
also hold for $\{I_j\}_{j=1}^n$ in a neighbourhood of an 
elliptic fixed point, under perturbations
depending also on an additional transverse component. In this paper we
study this situation in detail, taking particular interest in the possibility
of obtaining results which are uniform in the dimension of the transverse 
component. We shall consider perturbations
in which the additional transverse variable $\zeta=(\xi, \eta)\in\R^{2N}$, 
while the original phase space $\R^{2n}\ni z=(x,y)$ is a symplectic subspace
of the new enlarged phase space $\R^{2n}\times \R^{2N}$. We consider
Hamiltonians coupling $z$ and $\zeta$ of the form\footnote{It is 
possible to put the coupling between $z$ and $\zeta$ into either $f$
or $\Lambda$, and we make different choices depending upon which is
most convenient.}
$$H=H_0+\kappa\Lambda,\quad \Lambda=O(|\zeta|^2)\quad\hbox{as}\quad\zeta\to 0$$ 
and ask the question: under which conditions does Nekhoroshev 
exponential stability hold for $z\in\R^{2n}$ in a neighbourhood 
of the origin? (Notice that in contrast to some other discussions
we are not asking for the Nekhoroshev estimates
to hold for the full flow, only for the $z=\Pi_1(z,\zeta)$ part 
projected out of it.) The obvious perturbative problem arises 
by considering small $\kappa$; here our main theorem
\ref{gszsk} can be stated heuristically as:
\begin{quote}
{\em 
There exists a neighbourhood ${\cal N}$ of the origin in $\R^{2n}\times \R^{2N}$
and $\kappa_0$ such that for $0<\kappa<\kappa_0$ and initial
data in ${\cal N}$ exponential stability estimates like
\eqref{nee} hold for the flow projected onto $\R^{2n}$. All of
the neighbourhoods and estimates can be bounded explicitly and uniformly in $N$.
}
\end{quote}
(See also theorems \ref{agszsk} and \ref{gszskvar} 
for alternative formulations.) 
The proof of theorem \ref{gszsk} relies on a normal form lemma 
\ref{norm} which involves applying the method of averaging in a way which
couples $z$ and $\zeta$.

Counterintuitively perhaps, the case $\kappa\nearrow +\infty$ can sometimes
also be regarded as a perturbation of 
a Hamiltonian flow on $\R^{2n}$,  as we discuss in section \ref{scpot}. 
To be precise this is the case  when 
$\Lambda$ is such as  to
force the flow onto the $\R^{2n}\times\{0\}$
subspace for large $\kappa$, on which subspace  the dynamics is governed by 
the restricted Hamiltonian
$H_0(z)=H(z,0)$, that is,  motion in a constraining potential. 
(To ensure this, it is now required that $\Lambda=0$ 
if and only if $\zeta=0$; see section \ref{scpot} for the precise conditions).
In this case 
we have theorem \ref{gszk}, which can be stated heuristically as:
\begin{quote}
{\em 
Exponential stability estimates like
\eqref{nee} continue to hold for the  $\R^{2n}$ projection of
the flow in $\R^{2n}\times \R^{2N}$
determined by the Hamiltonian $H=H_0+\kappa\Lambda$,
with $\Lambda$ a constraining potential,
when $\kappa$ is sufficiently large (independent of $N$).
}
\end{quote}
In contrast to the small $\kappa$ case, these large $\kappa$ results 
in general come only with
the assurance that they hold for sufficiently large $\kappa$, but without
precise quantitative information on their domain of validity 
(at least in the absence of more special assumptions).
In fact theorem \ref{gszk} comes almost for free by combining general 
compactness results based on the Arzel\`{a}-Ascoli theorem
with the standard Nekhoroshev estimate for $H_0(z)$.
We provide details, (i) to explain the method in a simple case, 
(ii) to clarify how quantitative information on the domains can be
derived in special cases (theorem \ref{quantivers})
and (iii) to emphasize the difference with
the more involved case $\kappa\searrow 0$
which is treated in section \ref{nssk}. Finally we remark that stronger
results could be proved with the assumption that $\Lambda$ vanishes
faster than quadratically as $\zeta\to 0$, but we are not aware of
any likely applications in this case. 

Some results in a similar direction were obtained before
in \cite[p.~1713]{niederman}.
A crucial feature of our work is that we only do impose
conditions on the initial data rather then on the dynamics.
Moreover, although $N$ is finite,
all bounds are explicit and independent of $N$. This indicates that
it should be possible to obtain results at least for some
infinite  dimensional problems by the methods developed here, which will be the subject
of future work. Another possibility of generalizing our results
would consist of trying to relax the differentiability assumptions,
as in \cite{boune}.

%%%%%%%%%%%%%%%%%%%%%%%%%%%%%%%%%%%%%%%%%%%%%%%%%%%%%%%%%%%%%%%%%%%%%%%%%%%%%%%%%%%%%%%

\setcounter{equation}{0}

\section{Some notation} 

In general we will be concerned with real analytic Hamiltonians
$H=H(z, \zeta)$ depending on the variables $(z, \zeta)\in\R^{2n}\times\R^{2N}$
or $(z, \zeta)\in\C^{2n}\times\C^{2N}$. (By real analytic mapping from a complex
domain into another complex vector space,
we mean a complex analytic mapping which maps real vectors
into real vectors.)
Denoting $z=(z_1, \ldots, z_n)$ for $z_j=(x_j, y_j)\in\C^2$
for $1\le j\le n$ we write $I_j=(x_j^2+y_j^2)/2\in\C$, and also
$\zeta_j=(\xi_j, \eta_j)\in\C^2$ for $1\le j\le N$. Define the domains
\[ {\cal D}_{a,\,b,\,c}=\{(z, \zeta)\in\C^{2n}\times\C^{2N}:
   |I-{I^0}|<a, |z|<b, |\zeta|<c\} \]
for $a, b, c>0$, where ${I^0}\in\R^n$ is given and 
$$
|I-{I^0}|=\sum_{j=1}^n|I_j-{I_j^0}|,\quad
|z|^2=\sum_{j=1}^n(|x_j|^2+|y_j|^2)\quad\hbox{and}\quad
|\zeta|^2=\sum_{j=1}^N(|\xi_j|^2+|\eta_j|^2).
$$
The norm of a matrix $A\in\R^{n\times n}$ is the operator norm 
w.r.t. the $l_1$-norm $|I|=\sum_{j=1}^n |I_j|$. We always view $I$ 
as a function of $z$ and note the estimate, with 
$\tz_j=(\tx_j,\ty_j)$:  
\begin{eqnarray}\label{Iwz} 
   |I(\tz)-I(z)| & \le & \frac{1}{2}\,\Bigl(\sum_{j=1}^n|\tx_j-x_j|^2
   +|\ty_j-y_j|^2\Bigr)^{\frac{1}{2}}\Bigl(\sum_{j=1}^n|\tx_j+x_j|^2+
   |\ty_j+y_j|^2\Bigr)^{\frac{1}{2}}
   \nonumber
   \\ & \leq &\frac{1}{2}\,|\tz-z|\,(|\tz|^2+2|\tz||z|+|z|^2)^{\frac{1}{2}}
   =\frac{1}{{2}}\,|\tz-z|\,(|\tz|+|z|).
\end{eqnarray} 
The Hamiltonian vector field
generated by a function $f=f(z, \zeta)$ is written as $X_f$, and the
associated flow as $X^t_f$. We shall refer to integral curves of
$X_f$ also as integral curves of $f$ when no confusion seems likely.
The supremum norm of functions or vector fields on ${\cal D}_{a,\,b,\,c}$
is denoted by ${|\cdot|}_{a,\,b,\,c}$. For $\underline{r}=(r_1,r_2,r_3)$ we will
write ${\cal D}_{\underline{r}}={\cal D}_{\,r_1,\,r_2,\,r_3}$ and
${|\cdot|}_{\underline{r}}={|\cdot|}_{r_1,\,r_2,\,r_3}$.
Let $\Pi_1$ (resp. $\Pi_2$) be the orthogonal projection operator
onto the $\C^{2n}$ (resp. $\C^{2N}$) factor of $\C^{2n}\times\C^{2N}$.
We will refer to $\zeta=\Pi_2(z, \zeta)$ as the {\it transverse} component.
The symbols $C, C_1, C_2, \ldots$ are reserved for constants which are allowed
to depend only on $n$.

%%%%%%%%%%%%%%%%%%%%%%%%%%%%%%%%%%%%%%%%%%%%%%%%%%%%%%%%%%%%%%%%%%%%%%%%%%%%%%%%%%%%%%%%

\setcounter{equation}{0}

\section{One step improvement of the interaction term}

We start with an integrable Hamiltonian
$\langle\omega^0, I\rangle+\frac{1}{2}\,\langle A(I-{I^0}), I-{I^0}\rangle$
and a further Hamiltonian $\Lambda(\zeta)$
on $\R^{2N}$. We will introduce a coupling and use the following lemma
iteratively to successively reduce the interaction. In the proof we 
will sometimes abbreviate:
\begin{equation}\label{g0-def}
   h(z)=\langle\omega^0, I(z)\rangle\quad\hbox{and}\quad 
   g_0(z)=\frac{1}{2}\,\langle A(I(z)-{I^0}), I(z)-{I^0}\rangle.
\end{equation}

\begin{lemma}[Iteration step]\label{itstep-lem}
Consider the Hamiltonian
\[ H(z, \zeta)=\langle\omega^0, I(z)\rangle
   +\frac{1}{2}\,\langle A(I-{I^0}), I-{I^0}\rangle
   +g(z, \zeta)+f(z, \zeta)+\kappa\Lambda(\zeta), \] 
where $\omega^0, {I^0}\in\R^n$, $A\in\R^{n\times n}$
is a symmetric matrix and $T, \kappa>0$ are fixed
such that $T\omega^0\in 2\pi\Z^n$ holds.
The functions $g$ and $f$ are assumed to be real analytic on
an open set containing $\overline{{\cal D}}_{\underline{r}}$
for $\underline{r}=(r_1,r_2,r_3)$ with $r_1, r_2, r_3>0$,
whereas $\Lambda$
is assumed to be real analytic on an open set containing
$\{|\zeta|\leq r_3\}$. We suppose that
for some $\delta, \eps>0$ and some constant $C_\Lambda>0$,
\begin{enumerate}
\item ${|g|}_{\,\underline{r}}\le\delta$ and $\{g, h\}=0$,
\item ${|f|}_{\,\underline{r}}\le\eps$,
\item $|D\Lambda(\zeta)|\leq C_\Lambda |\zeta|$ for $|\zeta|\leq r_3$.
\end{enumerate}
If $\rho_1\in ]0, r_1[$, $\rho_2\in ]0, r_2[$, $\rho_3\in ]0, r_3[$
are such that
\begin{equation}\label{quanrel}
   \eps T<\frac{1}{9}\Big(\min\Big\{\frac{\rho_1}{r_2}, \rho_2, \rho_3\Big\}\Big)^2,
\end{equation}
then there exists a real analytic symplectic transformation
\[
\Phi:\,\,{\cal D}_{\underline{r}-\underline{\rho}}\to {\cal D}_{\underline{r}}
\]
such that, on ${\cal D}_{{\underline r}-\underline{\rho}}$,
\begin{equation}\label{HPhi-form}
   H\circ\Phi=\langle\omega^0, I(z)\rangle
+\frac{1}{2}\,\langle A(I-{I^0}), I-{I^0}\rangle
   +g_+(z, \zeta)+f_+(z, \zeta)+\kappa\Lambda(\zeta)
\end{equation}
and with the properties:
\begin{itemize}
\item[{\rm (a)}] \qquad
$\displaystyle {|\Phi-{\rm id}|}_{\,\underline{r}-\underline{\rho}}
\le\frac{3\,\eps T}{\min\{\frac{\rho_1}{r_2},\rho_2,\rho_3\}}$,
\item[{\rm (b)}] \qquad
${|{g_+}|}_{\,\underline{r}}\le\delta+\eps$ and $\{g_+, h\}=0$,
\item[{\rm (c)}] \qquad
$\displaystyle {|{f_+}|}_{\underline{r}-\underline{\rho}}
\le\bigg[\,\frac{6\|A\|r_1 r_2}{\rho_2}
+\frac{36\,(\delta+\eps)}{(\min\{\frac{\rho_1}{r_2},\rho_2,\rho_3\})^2}
+\frac{3\kappa C_{\Lambda}r_3}{2\rho_3}\bigg]\,\eps T.
$
\end{itemize}
\end{lemma}
{\bf Proof of lemma \ref{itstep-lem}\,} We start by averaging
over the flow generated by $h$: let
\begin{equation}\label{barf-def}
   \bar{f}(z, \zeta)=\frac{1}{T}\int_0^T (f\circ X_h^t)(z, \zeta)\,dt.
\end{equation}
Explicitly,
\begin{eqnarray*}
   X_h^t(z, \zeta) & = & (z_1(t),\ldots, z_n(t), \zeta),
   \\[1ex] z_j(t) & = & R_j(t)z_j,\quad z_j=(x_j, y_j),
   \\[1ex] R_j(t) & = & \left(\begin{array}{cc}
   \cos(\omega^0_j t) & \sin(\omega^0_j t) \\
   -\sin(\omega^0_j t) & \cos(\omega^0_j t) \\
   \end{array}\right).
\end{eqnarray*}
Since $T\omega^0\in 2\pi\Z^n$ we get $R_j(t+T)=R_j(t)$
and the flow $X_h^t$ is $T$-periodic.
In addition the matrices are real and $R_j^TR_j={\rm id}$, so that $|z(t)|=|z(0)|$ 
and $\I_j(t)=(x_j(t)^2+y_j(t)^2)/2=\I_j(0)$. 
Then $X_h^t$ leaves invariant
every domain ${\cal D}_{\underline{r}}$; in particular, $\bar{f}$
is well defined on ${\cal D}_{\underline{r}}$, the domain of $f$, and
${|\bar{f}|}_{\,\underline{r}}\leq {|f|}_{\,\underline{r}}\le\eps$. Now define
\beq\la{vaphi-def}
\varphi(z, \zeta)=\frac{1}{T}\int_0^T t\,
((f-\bar{f})\circ X_h^t)(z, \zeta)\,dt 
\eeq
which is well defined on ${\cal D}_{\,\underline{r}}$ and satisfies
\begin{equation}\label{varphi-prop}
   \{\varphi, h\}=f-\bar{f} \quad\hbox{and}\quad
   {|\varphi|}_{\,\underline{r}}\leq T{|f|}_{\,\underline{r}}\le\eps T.
\end{equation}
[To establish (\ref{varphi-prop}), we use
\begin{eqnarray*}
   \frac{d}{ds}\,\Big[s\,(f-\bar{f})\circ X_h^{t+s}\Big]
   & = & s\,\frac{d}{ds}\Big[(f-\bar{f})\circ X_h^{t+s}\Big]
   +(f-\bar{f})\circ X_h^{t+s}
   \\ & = & s\,\frac{d}{dt}\Big[(f-\bar{f})\circ X_h^{t+s}\Big]
   +(f-\bar{f})\circ X_h^{t+s},
\end{eqnarray*}
which upon integration $\int_0^T ds$ yields
\[ T(f-\bar{f})\circ X_h^{t+T}
   =\frac{d}{dt}\,\int_0^T s\,(f-\bar{f})\circ X_h^{t+s}\,ds
   +\int_0^T (f-\bar{f})\circ X_h^{t+s}\,ds. \]
Therefore
\begin{eqnarray*}
   \{\varphi, h\} & = & {\{\varphi, h\}\circ X_h^t\,\Big|}_{t=0}
   ={\frac{d}{dt}\Big(\varphi\circ X_h^t\Big)\,\Big|}_{t=0}
   =\frac{1}{T}\,\frac{d}{dt}\bigg(\int_0^T s\,(f-\bar{f})\circ X_h^{t+s}\,ds\bigg)\,\bigg|_{t=0}
   \\ & = & f-\bar{f}-\frac{1}{T}\int_0^T (f-\bar{f})\circ X_h^s\,ds
   \\ & = & f-\bar{f},
\end{eqnarray*}
(as a consequence of $X_h^T={\rm id}$ and the fact that $\bar{f}\circ X_h^s$
is independent of $s$ since
$$
\frac{d}{ds}\,\Big(\bar{f}\circ X_h^s\Big)=\frac{d}{ds}\,\bigg(
\frac{1}{T}\int_0^T f\circ X_h^{s+t}\,dt\bigg)
=\frac{1}{T}\int_0^T\frac{d}{dt}\,(f\circ X_h^{s+t})\,dt=0,
$$
so that the integral in the penultimate line is zero.)
The formula for $\varphi$ can be estimated in the obvious way
given the remarks already made on the action of $X^t_h$, completing
the proof of (\ref{varphi-prop}).]

Estimates for the derivatives of
$\varphi$ follow from Cauchy's theorem:
\begin{equation}\label{cauch}
{\bigg|\frac{\partial\varphi}{\partial z}\bigg|}_{\underline{r}-\underline{\rho}/3}
\leq \frac{3 {|\varphi|}_{\underline{r}}}{\min\{\frac{\rho_1}{r_2},\rho_2\}},\qquad
{\bigg|\frac{\partial\varphi}{\partial \zeta}\bigg|}_{\underline{r}-\underline{\rho}/3}
\leq \frac{3 {|\varphi|}_{\underline{r}}}{\rho_3},
\end{equation}
since $(z,\zeta)\in
{\cal D}_{\underline{r}-\underline{\rho}/3}$ and $|z-w|\leq\frac{1}{3}
\min\{\frac{\rho_1}{r_2}, \rho_2\}$ implies $(w,\zeta)\in 
{\cal D}_{\underline{r}}$. 
In fact, by (\ref{Iwz}), 
\begin{eqnarray*} 
   |I(w)-{I^0}| & \le & |I(w)-I(z)|+|I(z)-{I^0}|\le\frac{1}{2}\,|w-z|(|w-z|+2|z|)
   +r_1-\rho_1/3
   \\ & < & \frac{\rho_1}{6r_2}\Big(\rho_2/3+2(r_2-\rho_2/{3})\,\Big)+r_1-\rho_1/3<{r_1}. 
\end{eqnarray*}    
This implies bounds for the corresponding Hamiltonian vector field $X_\varphi$:
\begin{equation}\label{estdec}
{\big|\Pi_1 X_\varphi\big|}_{\underline{r}-\underline{\rho}/3}
\leq\frac{3\,\eps T}{\min\{\frac{\rho_1}{r_2},\rho_2\}},\qquad
{\big|\Pi_2 X_\varphi\big|}_{\underline{r}-\underline{\rho}/3}
\leq\frac{3\,\eps T}{\rho_3}.
\end{equation}

\begin{remark}
These Hamiltonian vector fields have, respectively, $2n$ and $2N$
components and the bounds (\ref{estdec}) hold using the Euclidean norm
with respect to these components; see \cite[Lemma 1]{lochak1} or
\cite[Prop.\ 3 in \S6]{nachbin} for an abstract treatment for maps
between Banach spaces.
\end{remark}
We now introduce
\begin{eqnarray}
   \Phi & = & X_\varphi^1,\qquad\qquad\hbox{(the time one map of the flow of $\varphi$)}
   \nonumber
   \\ g_+ & = & g+\bar{f},
   \label{g+-def}
   \\ f_+ & = & \int_0^1\{g_0+g+f_t, \varphi\}\circ X_\varphi^t\,dt\,
   +\,\kappa\,(\Lambda\circ\Phi-\Lambda),
   \label{f+-def}
\end{eqnarray}
where $g_0$ is as in (\ref{g0-def}) and
$f_t=tf+(1-t)\bar{f}$ for $t\in [0, 1]$.
To verify that (\ref{HPhi-form}) holds with the properties asserted,
observe that
\[ \frac{d}{dt}\Big[(g_0+g+f_t)\circ X_\varphi^t\Big]
   =\{g_0+g+f_t, \varphi\}\circ X_\varphi^t+(f-\bar{f})\circ X_\varphi^t, \]
and consequently
\begin{equation}\label{owsiat}
   (g_0+g+f)\circ\Phi-(g_0+g+\bar{f})
   =\int_0^1\{g_0+g+f_t, \varphi\}\circ X_\varphi^t\,dt
   +\int_0^1 (f-\bar{f})\circ X_\varphi^t\,dt.
\end{equation}
Since
\[ \frac{d}{dt}\Big(h\circ X_\varphi^t\Big)
   =-\,\{\varphi, h\}\circ X_\varphi^t=-\,(f-\bar{f})\circ X_\varphi^t
\]
by (\ref{varphi-prop}), it follows from (\ref{owsiat}) that
\[ (g_0+g+f+h)\circ\Phi-(g_0+g+\bar{f}+h)
    =\int_0^1\{g_0+g+f_t, \varphi\}\circ X_\varphi^t\,dt.
  \]
Thus
\begin{eqnarray*}
   H\circ\Phi & = & (h+g_0+g+f)\circ\Phi+\kappa\Lambda\circ\Phi
   \\ & = & h+g_0+g+\bar{f}+\kappa\Lambda+\kappa\,(\Lambda\circ\Phi-\Lambda)
   \\ & & +\,\int_0^1\{g_0+g+f_t, \varphi\}\circ X_\varphi^t\,dt
   \\ & = & h+g_0+g_+ +\kappa\Lambda+f_+
\end{eqnarray*}
which is the form of $H\circ\Phi$ asserted in (\ref{HPhi-form}),
with the functions $g_+$ and $f_+$ being defined
in (\ref{g+-def}) and (\ref{f+-def}), respectively.
To check the estimate (c) in the lemma we split up ${f_+}$ as follows:
\begin{eqnarray}
f_+ & = & \int_0^1\{g_0, \varphi\}\circ X_\varphi^t\,dt
+\int_0^1\{g+f_t, \varphi\}\circ X_\varphi^t\,dt
+\kappa\,(\Lambda\circ\Phi-\Lambda) \nonumber
\\ & = &f_{+,\,1}+f_{+,\,2}+f_{+,\,3}. \label{f+-defdec}
\end{eqnarray}
In order to derive the bounds for the $f_{+,\, j}$ quantities and to justify the
preceding calculation we summarize some mapping properties of the flows
in the following proposition,
thus also establishing statement (a) in the lemma since $\Phi=X_\varphi^1$.

\begin{pro}[{Mapping properties for the flows $X^t_\varphi$ and $X^t_{g_0}$}]
\label{mp}
Under the assumptions of lemma \ref{itstep-lem} the 
Hamiltonian flows generated by $\varphi$ and $g_0$ have the following
properties:
\begin{enumerate}

\item For real times $|t|\le 1$, the flow $X^t_\varphi$ satisfies
\begin{eqnarray}
   & & X_\varphi^t: {\cal D}_{\underline{r}-\underline{\rho}}
   \to {\cal D}_{\underline{r}-2\underline{\rho}/3},
   \label{mapprop1}
   \\[1ex] & & X_\varphi^t: {\cal D}_{\underline{r}-2\underline{\rho}/3}
   \to {\cal D}_{\underline{r}-\underline{\rho}/3}\quad\hbox{and}
   \label{mapprop2}
   \\[1ex] & &  {|X_\varphi^t-{\rm id}|}_{\,\underline{r}-\underline{2\rho}/3}
   \le {|X_\varphi|}_{\,\underline{r}-\underline{\rho}/3}\,|t|
   \leq\frac{3\,\eps T}{\min\{\frac{\rho_1}{r_2},\rho_2,\rho_3\}},
   \label{ftfp}
\end{eqnarray}
and for complex times $t$ such that
\begin{equation}\label{lamdef}
   |t|<\lambda\quad\mbox{for}\quad\lambda
   =\frac{1}{18\,\eps T}\biggl(\min\Big\{\frac{\rho_1}{r_2}, \rho_2, \rho_3\Big\}\biggr)^2
\end{equation}
the flow $X^t_{\varphi}$ is analytic on ${\cal D}_{\underline{r}-\underline{\rho}/6}$
and satisfies
\begin{eqnarray}
   & & X_\varphi^t: {\cal D}_{\underline{r}-2\underline{\rho}/3}
   \to {\cal D}_{\underline{r}-\underline{\rho}/2}\subset
   {\cal D}_{\underline{r}-\underline{\rho}/3}.
   \label{suf1}
\end{eqnarray}

\item For complex times $t$ such that
\begin{equation}\label{taudef}
   |t|<\tau\quad\mbox{for}\quad\tau=\frac{\rho_2}{3\|A\|\,r_1 r_2}
\end{equation}
the flow $X^t_{g_0}$ is analytic on ${\cal D}_{\underline{r}-\underline{\rho}/3}$
and satisfies
\begin{eqnarray}
   & & X_{g_0}^t: {\cal D}_{\underline{r}-2\underline{\rho}/3}
   \to {\cal D}_{\underline{r}-\underline{\rho}/3},
   \label{klbgr1}
   \\[1ex] & & X_{g_0}^t: {\cal D}_{\underline{r}-\underline{\rho}/3}
   \to {\cal D}_{\underline{r}}\,.
   \label{klbgr2}
\end{eqnarray}
\end{enumerate}
\end{pro}
{\bf Proof of Proposition \ref{mp}\,}
(i) To begin with, the equation $(d/dt)X_\varphi^t=X_\varphi(X_\varphi^t)$ reads
\[ (\dot{z}(t), \dot{\zeta}(t))=(\Pi_1 X_\varphi,\Pi_2X_\varphi)(z(t),\zeta(t)) \]
where $(z(t), \zeta(t))=X_\varphi^t(z(0), \zeta(0))$
for some fixed $(z(0), \zeta(0))\in {\cal D}_{\underline{r}-\underline{\rho}}$. This
implies, by (\ref{estdec}), that
\begin{equation}\label{ablest}
|\dot{z}(t)|\le\frac{3\,\eps T}{\min\{\frac{\rho_1}{r_2},\rho_2\}}
\quad\hbox{and}\quad
|\dot{\zeta}(t)|\leq\frac{3\,\eps T}{\rho_3},
\end{equation}
at least as long as the solution stays in ${\cal D}_{\underline{r}-\underline{\rho}/3}$,
during which time
\[ |\dot{I}|=\Big|\sum_{j=1}^n (x_j\dot{x}_j+y_j\dot{y}_j)\Big|
   \le |z(t)||\dot z(t)|
\leq\frac{3\,r_2\eps T}{\min\{\frac{\rho_1}{r_2}, \rho_2\}}. \]
Writing $\I(t)=I(z(t))$ with $z(0)\in {\cal D}_{\underline{r}-\underline{\rho}}$
we deduce from (\ref{quanrel}) that for $|t|\leq 1$
\[ |\I(t)-{I^0}|\le |\I(t)-\I(0)|+|\I(0)-{I^0}|\le
\frac{3\,r_2|t|\eps T}{\min\{\frac{\rho_1}{r_2}, \rho_2\}}
+r_1-\rho_1<r_1-\frac{2}{3}\,\rho_1. \]
Furthermore, using (\ref{quanrel}) again,
\[ |z(t)|\le |z(0)|+\frac{3\,|t|\eps T}{\min\{\frac{\rho_1}{r_2}, \rho_2\}}
   \le r_2-\rho_2+\frac{3\,|t|\eps T}{\min\{\frac{\rho_1}{r_2}, \rho_2\}}
   <r_2-\frac{2}{3}\,\rho_2 \]
and
\[ |\zeta(t)|\le |\zeta(0)|+\frac{3\,|t|\eps T}{\rho_3}\le r_3-\rho_3
   +\frac{3\,|t|\eps T}{\rho_3}
   <r_3-\frac{2}{3}\,\rho_3. \]
This argument shows in particular that if the $\rho_j$ are chosen
in accordance with (\ref{quanrel}), then the solution starting in
${\cal D}_{\underline{r}-\underline{\rho}}$ will remain
in ${\cal D}_{\underline{r}-2\underline{\rho}/3}$ for all times $|t|\le 1$.
This proves (\ref{mapprop1}), and verification of (\ref{mapprop2}) is analogous.
Moreover, (\ref{ftfp}) follows from (\ref{ablest}).

For the complex case, since $X_{\varphi}$ is analytic,
the flow $(X_{\varphi}^t)$ is defined locally and is locally analytic
on $\C^{2n}\times\C^{2N}$ and for complex $t$. To find for which
$t\in\C$ and between which domains this is true, we
just repeat the argument that led
to (\ref{mapprop1}), and it is found that for $|t|<\lambda$
with $\lambda$ as in (\ref{lamdef}) the flow is well defined,
analytic and satisfies (\ref{suf1}).

(ii) Again, since $X_{g_0}$ is analytic,
the flow $(X_{g_0}^t)$ is defined locally,
and is locally analytic,
on $\C^{2n}\times\C^{2N}$  for complex $t$. Observe that
\[ \frac{d}{dt}\Big(I\circ X_{g_0}^t\Big)
   =\{I, g_0\}\circ X_{g_0}^t=0 \]
for the function $I=I(z)$, since $g_0=g_0(I)$ only depends on $z$ through
$I=(I_1,\dots I_n)$.
In addition, since $g_0$ is independent of
$\zeta=\Pi_2(z, \zeta)$ we have
\[ \frac{d}{dt}\Big(\zeta\circ X_{g_0}^t\Big)
   =\{\zeta, g_0\}\circ X_{g_0}^t=0. \]
In other words, both $I$ and $\zeta$
are preserved by the flow, so that restrictions on the time which
ensure (\ref{klbgr1})-(\ref{klbgr2}) arise only  from the condition
on $z$. To prove (\ref{klbgr2})
for instance, write $(z(t), \zeta(t))=X_{g_0}^t(z(0), \zeta(0))$
and  $\I(t)=I(z(t))$. Then by the foregoing observation:
\begin{eqnarray*}
   |\I(t)-{I^0}| & = & |\I(0)-{I^0}|<r_1-\rho_1/3<r_1,
   \\ |\zeta(t)| & = & |\zeta(0)|<r_3-\rho_3/3<r_3,
\end{eqnarray*}
for all times, provided that initially $(z(0), \zeta(0))
\in {\cal D}_{\underline{r}-\underline{\rho}/3}$.
Furthermore,
\[ |\dot{z}(t)|\le\Big|\frac{d}{dt}\,X_{g_0}^t\Big|
   =|X_{g_0}(X_{g_0}^t)|\le {|X_{g_0}|}_{\,\underline{r}} \]
as long as the flow stays in ${\cal D}_{\underline{r}}$.
Using the definition of $g_0$ in (\ref{g0-def}),
for any unit $2n$ vector $(\mathbf{a}, \mathbf{b})=(a_1,\dots,a_n,b_1,\dots b_n)$ 
we can estimate
\begin{eqnarray*} 
|(\mathbf{a}\cdot\nabla_x+\mathbf{b}\cdot\nabla_y){g_0}|
& = & |\sum_{i, j=1}^n A_{ij}(I-{I^0})_i(a_jx_j+b_jy_j)|
\\ & \le & \|A\|\,|I-{I^0}|\,\max_{1\le j\le n} |a_jx_j+b_jy_j|\le\|A\|\,|I-{I^0}|\,|z|.
\end{eqnarray*} 
So, maximizing over the unit vector, we can bound
the Euclidean norm for $X_{g_0}$ as
\[ |X_{g_0}(z, \zeta)|\le\|A\|\,|I-{I^0}|\,|z| \]
(using the $l^1$ operator norm on $A$).
It follows that
\[ {|X_{g_0}|}_{\,\underline{r}}\le\|A\|\,r_1 r_2
   \quad\mbox{for\,\,any}\quad\underline{r}. \]
Hence the desired bound $|z(t)|<r_2$ is obtained by inserting
(\ref{taudef}) into the estimate:
\[ |z(t)|\le |z(0)|+|t|\,{|X_{g_0}|}_{\,\underline{r}}
   <r_2-\frac{1}{3}\,\rho_2+\|A\|\,r_1 r_2\,|t|<r_2. \]
To summarize, it has been shown that (\ref{klbgr2}) is verified for $|t|<\tau$,
and (\ref{klbgr1}) follows in the same way.
\qed\bigskip

\smallskip

\noindent
{\bf Continuation of proof of lemma \ref{itstep-lem}\,}
So far the statements (\ref{HPhi-form}) and (a) of the lemma are proved.
Next, notice that the first assertion in (b) follows immediately from
the definition of $g_+$ in (\ref{g+-def}), and
the assumption ${|f|}_{\,\underline{r}}\le\eps$.
To establish the second assertion in (b)
we need to prove that $\{\bar{f}, h\}=0$
(in view of $g_+=g+\bar{f}$ and $\{g, h\}=0$)
which follows directly from the definition (\ref{barf-def}):
\[
\{\bar{f}, h\}=\frac{1}{T}\int_0^T \{f\circ X_h^t, h\}\,dt
=\frac{1}{T}\int_0^T\frac{d}{dt}\,(f\circ X_h^t)\,dt =0.
\]
To complete the proof of the lemma it remains to verify (c),
which is now done by estimating each of
the three terms in (\ref{f+-defdec}).

{\it Estimation of $f_{+,\,1}$:}
As a consequence of (ii) in the previous proposition, the function
\[ F(t)=\varphi\circ X_{g_0}^t(z, \zeta) \]
is analytic for complex times $t$ as in (\ref{taudef})
and for $(z,\zeta)\in {\cal D}_{\underline{r}-\underline{\rho}/3}$,
since $\varphi$ is defined on ${\cal D}_{\,\underline{r}}$.
Then by Cauchy's estimate
\[ |\{{g_0}, \varphi\}(z, \zeta)|=|F'(0)|
   \le\frac{2}{\tau}\,\sup_{|t|=\tau/2}\,|F(t)| \]
for every $(z, \zeta)\in {\cal D}_{\underline{r}-2\underline{\rho}/3}$.
To bound $F(t)=\varphi\circ X_{g_0}^t$ we just observe that by (\ref{klbgr1})
and (\ref{varphi-prop}),
\[ {|\varphi\circ X_{g_0}^t|}_{\,\underline{r}-2\underline{\rho}/3}
   \le {|\varphi|}_{\,\underline{r}-\underline{\rho}/3}
   \le {|\varphi|}_{\,\underline{r}}\le\eps T, \]
which leads to the estimate
${|\{{g_0},\varphi\}|}_{\,\underline{r}-2\underline{\rho}/3}\leq 2\eps T/\tau$.
Hence, by (\ref{mapprop1}) in the previous proposition,
\[ {|f_{+,\,1}|}_{\,\underline{r}-\underline{\rho}}
   =\bigg|\int_0^1\{g_0,\varphi\}
   \circ X_\varphi^t\,dt\bigg|_{\,\underline{r}-\underline{\rho}}
   \le {|\{{g_0},\varphi\}|}_{\,\underline{r}-2\underline{\rho}/3}
   \le\frac{6\|A\|\,r_1 r_2\eps T}{\rho_2}. \]

{\it Estimation of $f_{+,\,2}$:}
Next, to bound $f_{+,\,2}=\int_0^1\{g+f_t, \varphi\}\circ X_\varphi^t\,dt$
we proceed in a similar fashion, but using the flow $X^t_\varphi$ in place
of $X^t_{g_0}$. To treat the first term in the integral define
\[ G(t)=g\circ X_\varphi^t(z, \zeta) \]
where $(z, \zeta)\in {\cal D}_{\underline{r}-2\underline{\rho}/3}$ is fixed.
By (i) in the previous proposition this is analytic for complex times $t$
as in (\ref{lamdef}), so that Cauchy's estimate gives
\[ {|\{g, \varphi\}|}_{\,\underline{r}-2\underline{\rho}/3}=|G'(0)|
   \le\frac{2}{\lambda}\,\sup_{|t|=\lambda/2}\,
   {|G(t)|}_{\,\underline{r}-2\underline{\rho}/3}. \]
By (\ref{suf1}), $G$ is bounded as
${|G(t)|}_{\,\underline{r}-2\underline{\rho}/3}
\leq {|g|}_{\,\underline{r}-\underline{\rho}/3}\le\delta$
for these $t$, leading to the overall bound
\[ \Big|\int_0^1\{g, \varphi\}\circ X_\varphi^t\,dt\Big|_{\,\underline{r}-\underline{\rho}}
   \le {|\{g, \varphi\}|}_{\,\underline{r}-2\underline{\rho}/3}\le 2\delta/\lambda. \] 
The second term in the integral defining $f_{+,\,2}$ is handled in exactly
the same way, leading to the same bound with $\delta$ replaced by $\eps$,
since ${|f_t|}_{\,\underline{r}}\le\eps$ for $t\in [0, 1]$.
Therefore altogether
$$ {|f_{+,\,2}|}_{\,\underline{r}-\underline{\rho}}
   \le\frac{2(\delta+\eps)}{\lambda}
   =\frac{36\,\eps T(\delta+\eps)}
   {{(\min\{\frac{\rho_1}{r_2}, \rho_2, \rho_3\})}^2}.
$$

{\it Estimation of $f_{+,\,3}$:}
The last contribution to $f_+$ arises from
$f_{+,\,3}=\kappa\,(\Lambda\circ\Phi-\Lambda)$.
By definition of $\Phi=X_\varphi^1$ this can be rewritten as
\beq\la{ef3}
   f_{+,\,3}=\kappa\,(\Lambda\circ X_\varphi^1-\Lambda\circ X_\varphi^0)
   =\kappa\int_0^1\frac{d}{dt}\,(\Lambda\circ X_\varphi^t)\,dt
   =\kappa\int_0^1 \{\Lambda, \varphi\}\circ X_\varphi^t\,dt,
\eeq
so that, using (\ref{mapprop1}) and $\{\Lambda,\varphi\}=\langle D\Lambda, \Pi_2 X_\varphi\rangle$
(the latter due to $\Lambda=\Lambda(\zeta)$), we deduce
$$ {|f_{+,\,3}|}_{\,\underline{r}-\underline{\rho}}
   \le \kappa\,{|\langle D\Lambda, \Pi_2 X_\varphi\rangle|}_{\,\underline{r}-2\underline{\rho}/3}
   \le\kappa\,C_{\Lambda}r_3\,{|\Pi_2 X_\varphi|}_{\,\underline{r}-2\underline{\rho}/3}.
$$
Since only the $\zeta_j=(\xi_j, \eta_j)$ derivatives of $\varphi$ contribute to $\Pi_2 X_\varphi$,
this can be combined with Cauchy's estimate as
\[
    {|f_{+,\,3}|}_{\,\underline{r}-\underline{\rho}}
    \le \frac{3\kappa C_{\Lambda}r_3}{2\rho_3}\,{|\varphi|}_{\,\underline{r}}
    \le\frac{3\kappa C_{\Lambda}r_3}{2\rho_3}\,\eps T
\]
by (\ref{varphi-prop}). If we add together these bounds
on ${|f_{+,\,j}|}_{\,\underline{r}-\underline{\rho}}$,
then (c) is obtained. \qed\bigskip

%%%%%%%%%%%%%%%%%%%%%%%%%%%%%%%%%%%%%%%%%%%%%%%%%%%%%%%%%%%%%%%%%%%%%%%%%%%%%%%%%%%%%%%

\setcounter{equation}{0}

\section{Transformation to normal form}

We iterate Lemma \ref{itstep-lem} $m$ times to prove the following result.

\begin{lemma}[Normal form]\label{norm}
Consider the Hamiltonian
\[ H(z, \zeta)=\langle\omega^0, I\rangle+\frac{1}{2}\,\langle A(I-{I^0}), I-{I^0}\rangle
   +g(z, \zeta)+f(z, \zeta)+\kappa\Lambda(\zeta), \]
where $\omega^0, {I^0}\in\R^n$, $A\in\R^{n\times n}$ is a
symmetric matrix and $T, \kappa>0$ are fixed
such that $T\omega^0\in 2\pi\Z^n$ holds.
The functions $g$ and $f$ are assumed to be real analytic
on an open neighbourhood of
$\overline{\cal D}_{3\underline{r}}$, and $\Lambda$
is assumed to be real analytic on an open neighbourhood of
$\{|\zeta|\leq 3r_3\}$. We suppose that
\begin{enumerate}
\item ${|g|}_{\,3\underline{r}}\le\delta$ and $\{g, h\}=0$,
\item ${|f|}_{\,3\underline{r}}\le\eps$,
\item $|D\Lambda(\zeta)|\le C_\Lambda |\zeta|$ for $|\zeta|\le 3r_3$,
\item $r_1<2r_2^2$ and $r_1<2r_2r_3$,
\item $\displaystyle m^2\eps T<\frac{r_1^2}{81\,r_2^2}$,
\item $\displaystyle 54 m\|A\|r_1T+\frac{324\,(\delta+2\eps)m^2r_2^2T}{r_1^2}
+\frac{9\kappa C_{\Lambda}mT}{2}\leq\frac{1}{2}$
\end{enumerate}
for some $\delta, \eps>0$ and $C_\Lambda>0$.
Then there exists a real analytic symplectic transformation
\[
\Psi: {\cal D}_{\,2\underline{r}}\to {\cal D}_{\,3\underline{r}}
\]
such that, on ${\cal D}_{{2\underline r}}$,
\[ H\circ\Psi=\langle\omega^0, I\rangle+\frac{1}{2}\,\langle A(I-{I^0}),
   I-{I^0}\rangle+\hat g(z, \zeta)+\hat f(z, \zeta)+\kappa\Lambda(\zeta) \]
and with the properties:
\begin{itemize}
\item[{\rm (a)}] \qquad $\displaystyle {|\Psi-{\rm id}|}_{\,2\underline{r}}
\le\frac{18\,m r_2}{r_1}\,\eps T$,
\item[{\rm (b)}] \qquad ${|{\hat g}|}_{\,2\underline{r}}
\le\delta+2\eps$ and $\{\hat g, h\}=0$,
\item[{\rm (c)}] \qquad ${|{\hat f}|}_{\,2\underline{r}}\le  2^{-m}\eps$.
\end{itemize}
\end{lemma}
{\bf Proof\,} We apply the iterative lemma (Lemma \ref{itstep-lem})
$m$ times, where at the $j^{th}$ stage $\underline{r}$ is taken to be
$3\underline{r}-j\underline{r}/m$ and
$\underline{\rho}=\underline{r}/m$ with $j=0, \ldots, m-1$.
For $j=0$ we need to check (\ref{quanrel}), which reads as
\[ \eps T<\frac{1}{9\,m^2}\Big(\min\Big\{\frac{r_1}{3r_2}, r_2, r_3\Big\}\Big)^2. \]
According to (iv) we have $\min\{\frac{r_1}{3r_2}, r_2, r_3\}=\frac{r_1}{3r_2}$
and the condition becomes
\[ \eps T<\frac{1}{81\,m^2}\frac{r_1^2}{r_2^2} \]
which is verified by (v). Thus Lemma \ref{itstep-lem} yields
a real analytic symplectic transformation
$\Phi_1:\,\,{\cal D}_{3\underline{r}-\underline{r}/m}\to {\cal D}_{3\underline{r}}$
such that, on ${\cal D}_{3{\underline r}-\underline{r}/m}$,
\[ H\circ\Phi_1=\langle\omega^0, I\rangle+\frac{1}{2}\,\langle A(I-{I^0}), I-{I^0}\rangle
   +g_1(z, \zeta)+f_1(z, \zeta)+\kappa\Lambda(\zeta) \]
and moreover:
\begin{itemize}
\item\qquad
$\displaystyle {|\Phi_1-{\rm id}|}_{\,3\underline{r}-\underline{r}/m}
\le\frac{9\,m r_2}{r_1}\,\eps T$,
\item\qquad
${|{g_1}|}_{\,3\underline{r}}\le\delta+\eps$ and $\{g_1, h\}=0$,
\item\qquad
$\displaystyle {|{f_1}|}_{3\underline{r}-\underline{r}/m}
\le\bigg[\,54m\|A\| r_1
+\frac{324\,(\delta+\eps)m^2 r_2^2}{r_1^2}
+\frac{9\kappa C_{\Lambda}m}{2}\bigg]\,\eps T\le\frac{\eps}{2},
$
\end{itemize}
the latter in view of (vi). Put $\Psi_1=\Phi_1$.
For the induction step assume that we have constructed
a real analytic symplectic transformation $\Psi_j$ such that,
on ${\cal D}_{\,3\underline{r}-j\underline{r}/m}$,
$$
   H\circ\Psi_j=\langle\omega^0, I\rangle+\frac{1}{2}\,\langle A(I-{I^0}), I-{I^0}\rangle
   +g_j(z, \zeta)+f_j(z, \zeta)+\kappa\Lambda(\zeta)
$$
with
\begin{itemize}
\item\qquad $\displaystyle {|\Psi_j-{\rm id}|}_{\,3\underline{r}-j\underline{r}/m}
\le\frac{9\,m{r_2}}{r_1}\,\eps T\,\sum_{i=0}^{j-1} 2^{-i}$,
\item\qquad $\displaystyle {|{g_{j}}|}_{\,3\underline{r}-(j-1)\underline{r}/m}
\le\delta+\eps\sum_{i=0}^{j-1} 2^{-i}$ and $\{g_j, h\}=0$,
\item\qquad $ {|{f_{j}}|}_{\,3\underline{r}-j\underline{r}/m}\le  2^{-j}\eps$.
\end{itemize}
In order to apply the iterative lemma to this Hamiltonian 
(and with $\eps$ replaced by
$2^{-j}\eps$ and $\delta$ replaced by $\delta+\eps\sum_{i=0}^{j-1} 2^{-i}$),
we have to see that (\ref{quanrel}) holds, which reads as
\begin{equation}\label{applcon}
   2^{-j}\eps T<\frac{1}{9\,m^2}\Big(\min\Big\{\frac{r_1}{r_2(3-j/m)}, r_2, r_3\Big\}\Big)^2.
\end{equation}
Since $\frac{r_1}{r_2(3-j/m)}\le\frac{r_1}{2r_2}\le\min\{r_2, r_3\}$ by (iv),
(\ref{applcon}) reduces to
\[ 2^{-j}\eps T<\frac{1}{9\,m^2}\Big(\frac{r_1}{r_2(3-j/m)}\Big)^2\,,\]
which is a consequence of (v). Therefore Lemma \ref{itstep-lem} applies,
yielding a real analytic symplectic transformation
\[
\Phi_{j+1}: {\cal D}_{\,3\underline{r}-(j+1)\underline{r}/m}\to
{\cal D}_{\,3\underline{r}-j\underline{r}/m}
\]
such that, on ${\cal D}_{\,3\underline{r}-(j+1)\underline{r}/m}$,
\[ H\circ\Psi_j\circ\Phi_{j+1}=\langle\omega^0, I\rangle+\frac{1}{2}\,\langle A(I-{I^0}), I-{I^0}\rangle
   +g_{j+1}(z, \zeta)+f_{j+1}(z, \zeta)+\kappa\Lambda(\zeta) \]
and furthermore by the hypotheses:
\begin{itemize}
\item\qquad $\displaystyle {|\Phi_{j+1}-{\rm id}|}_{\,3\underline{r}-(j+1)\underline{r}/m}
\le\frac{3\cdot 2^{-j}\eps Tm}{(\frac{r_1}{r_2(3-j/m)})}
\le\frac{9\,mr_2}{r_1}\,2^{-j}\eps T$,
\item\qquad $\displaystyle {|{g_{j+1}}|}_{\,3\underline{r}-j\underline{r}/m}
\leq\delta+\eps\sum_{i=0}^{j-1} 2^{-i}+2^{-j}\eps=\delta+\eps\sum_{i=0}^j 2^{-i}$
and $\{g_{j+1}, h\}=0$,
\item\qquad $\displaystyle {|{f_{j+1}}|}_{\,3\underline{r}-(j+1)\underline{r}/m}
\leq\bigg[\,6m\|A\| r_1(3-j/m)^2
+\frac{36\,(\delta+\eps \sum_{i=0}^j 2^{-i})m^2}{(\frac{r_1}{r_2(3-j/m)})^2}
+\frac{3\kappa C_{\Lambda}(3-j/m)m}{2}\bigg]\,2^{-j}\eps T
\\ {} \hspace{9em} \le 2^{-(j+1)}\eps$.
\end{itemize}
Now define $\Psi_{j+1}=\Psi_j\circ\Phi_{j+1}$ and estimate
\begin{eqnarray*}
   {|\Psi_{j+1}-{\rm id}|}_{\,3\underline{r}-(j+1)\underline{r}/m}
   & \le & {|(\Psi_j-{\rm id})\circ\Phi_{j+1}|}_{\,3\underline{r}-(j+1)\underline{r}/m}
   +{|\Phi_{j+1}-{\rm id}|}_{\,3\underline{r}-(j+1)\underline{r}/m}
   \\ & \le & {|\Psi_j-{\rm id}|}_{\,3\underline{r}-j\underline{r}/m}
   +\frac{9\,mr_2}{r_1}\,2^{-j}\eps T
   \le\frac{9\,m{r_2}}{r_1}\,\eps T\,\sum_{i=0}^j 2^{-i}
\end{eqnarray*}
to deduce that the inductive assumptions hold also at this step.
The process terminates at $j=m-1$ and we can define $\hat{g}=g_{m-1}$,
$\hat{f}=f_{m-1}$, and $\Psi=\Psi_{m-1}$. \qed\bigskip

%%%%%%%%%%%%%%%%%%%%%%%%%%%%%%%%%%%%%%%%%%%%%%%%%%%%%%%%%%%%%%%%%%%%%%%%%%%%%%%%%%%%%%%

\setcounter{equation}{0}

\section{Nekhoroshev stability in the case $N=0$}

We recall the statement of Nekhoroshev stability in the case $N=0$, 
so that only the $z$ component appears. 
We assume that the initial values $z(0)=(x(0),y(0))\in\R^{2n}$ 
are close to the equilibrium point
$(0, 0)$ for the real analytic Hamiltonian:
\[  H(z)=\langle\alpha, I(z)\rangle+\frac{1}{2}\,\langle AI(z), I(z)\rangle+f(z),
    \quad\hbox{with}\quad\langle AI, I\rangle\geq\frac{1}{M}\,|I|^2
    =\frac{1}{M}\Big(\sum_{j=1}^n |I_j|\Big)^2 \] 
and $f(z)={\cal O}(z^5)$ for $|z|\to 0$. 
In that case we have the following theorem:

\begin{theorem}\label{classicalthm} 
There exist positive numbers $K, k, a$ (depending on $n$, $\alpha$, $M$ 
and $\|A\|$) and $\theta_0$ (depending on $n$, $\alpha$, $M$, $\|A\|$ and $f$) 
with the following properties. If $\I(0)=I(z(0))$ is such that 
$|\I(0)|=\theta^2$ for some $0<\theta\le\theta_0$, then 
$\I(t)=I(z(t))$ satisfies
\[ |\I(t)-\I(0)|\leq K\theta^{2+a}
   \quad\hbox{for}\quad |t|\le e^{\frac{k}{\theta^a}}. \] 
\end{theorem}
\noindent
{\bf Proof\,} This is the classical Nekhoroshev bound of \cite{nek77}
for the case of an elliptic equilibrium, see \cite{gfb,lochak1,nied98,poschel}.
The proof can also be extracted from the proof of our theorem \ref{gszk} below,
although not in its putative sharpest form (with 
$a=\frac{1}{2n})$.
\qed\bigskip

%%%%%%%%%%%%%%%%%%%%%%%%%%%%%%%%%%%%%%%%%%%%%%%%%%%%%%%%%%%%%%%%%%%%%%%%%%%%%%%%%%%%%%%%

\setcounter{equation}{0}

\section{Constrained motion: the case of large $\kappa$}
\la{scpot}

In this section we consider the case in which there is a 
transverse variable, $\zeta\in\R^{2N}$,
which is subject to a strong constraining potential. Precisely, we consider
real analytic Hamiltonians of the form
\begin{equation}\label{ham}
   H(z,\zeta)=H_0(z)+\kappa\Lambda(z, \zeta)\,,
\end{equation}
in the limit $\kappa\to +\infty$, assuming that
\begin{equation}\label{coer0}
\Lambda\geq 0\qquad\hbox{and }\quad
\Lambda(z, \zeta)=0\;\hbox{ if and only if}\;\zeta=0.
\end{equation}
The idea is that in this limit $\Lambda$  forces 
the motion onto the set $\Lambda=0$, thus dynamically enforcing the
constraint $\zeta=0$. We will work under the assumption
that there exist positive numbers $c_0,c_1,p$ such that for real $(z,\zeta)$
\begin{equation}\label{coer1}
   H(z,\zeta)\geq c_0 (|z|^p+\kappa|\zeta|^p)-c_1\,,
\end{equation}
and also that for all $R>0$ there exists $c_2(R)>0$ such that
\begin{equation}\label{coer2}
   \left|\frac{\partial \Lambda}{\partial z}(z,\zeta)\right |
   \leq c_2(R)\Lambda (z,\zeta)\quad\hbox{for}\quad |z|\leq R\,.
\end{equation}

\begin{remark}
The numbers $c_0,c_1,c_2,p$, and hence the bound \eqref{coer2},
are assumed to be independent of $\kappa$.
\end{remark}

The first result does not require analyticity:

\begin{lemma}\label{constr}
Assume that $H$ is a $C^{1,1}$ function (i.e.~$C^1$ with Lipschitz
derivative) of the form \eqref{ham}, also
verifying \eqref{coer0}-\eqref{coer2}.
Let there be given real initial data $(z^\kappa(0),\zeta^\kappa(0))$ 
for $\kappa\geq\kappa_0$, such that
\begin{enumerate}
\item
$z^\kappa(0)\to z(0)$ as $\kappa\to +\infty$;
\item
$\kappa\Lambda(z^\kappa(0),\zeta^\kappa(0))\to 0$ as $\kappa\to +\infty$;
\item
$\sup_{\kappa\geq\kappa_0} H(z^\kappa(0),\zeta^\kappa(0))=E<\infty$.
\end{enumerate}
Then there exist, for each $\kappa\geq\kappa_0$,
global integral curves $(z^\kappa(t),\zeta^\kappa(t))$
of $H$ which have the property that
\begin{equation}\la{urth0}
\lim_{\kappa\to +\infty}\max_{|t|\leq T}\bigl(|z^\kappa(t)-z(t)|+|\zeta^\kappa(t)|\bigr)=0
\end{equation}
for any $T>0$, where $z(t)$ is an integral curve of the Hamiltonian $H(z, 0)=H_0(z)$.
Furthermore
\begin{equation}\label{urth}
   \lim_{\kappa\to +\infty}\max_{|t|\leq T}\,
   \Bigl(|{H_0}(z^\kappa(t))-{H_0}(z^\kappa(0))|
   +\kappa\,\Lambda (z^\kappa(t),\zeta^\kappa(t))\Bigr)=0\,.
\end{equation}
\end{lemma}
{\bf Proof of lemma \ref{constr}\,}
The co-ercivity in (\ref{coer1}) together with energy conservation
implies the bound
\begin{equation}\label{zzet-bd}
   \kappa|\zeta^\kappa(t)|^p+|z^\kappa(t)|^p\leq \frac{c_1+E}{c_0}\,,
\end{equation} 
which is uniform in $\kappa$, and shows that $\zeta^\kappa(t)
=O({{\kappa}^{-\frac{1}{p}}})$,
uniformly in $t$. To obtain compactness for $z^\kappa(t)$ we use the $z$ 
component of the differential equation, i.e.
$$
\frac{d}{dt}\,z^\kappa=\Pi_1\, X_{H}(z^\kappa,\zeta^\kappa)\,,
$$
conservation of energy, \eqref{coer2} and (\ref{zzet-bd}) 
to deduce that $\dot z^\kappa(t)$ is bounded, uniformly in $t$ and
$\kappa\geq\kappa_0$. It follows from the Arzel\`{a}-Ascoli theorem that there
exists a subsequence converging uniformly on bounded intervals
$[-T, T]$ to a continuous limit $z=z(t)$. To prove that this limit
is an integral curve of $H_0$ we consider the integrated
form of the equation:
\beq\la{intf}
z^\kappa(t)=z^\kappa(0)
+\int_0^t\, \Pi_1\,X_{H}(z^\kappa(s),\zeta^\kappa(s))\,ds
=\int_0^t\, \,\Bigl[X_{H_0}(z^\kappa(s))
+\kappa\,\Pi_1\,X_{\Lambda}(z^\kappa(s),\zeta^\kappa(s))\Bigr]
\,ds
\,.
\eeq
Notice first that it is possible to
take the limit of this equation once we know \eqref{urth} holds, on
account of \eqref{coer2}. So we first prove \eqref{urth}.
Energy conservation 
$H_0(z^\kappa(t))+\kappa\Lambda(z^\kappa(t), \zeta^\kappa(t))
   =H_0(z^\kappa(0))+\kappa\Lambda(z^\kappa(0), \zeta^\kappa(0))
$
and the assumptions (i), (ii) imply that 
$\lim_{\kappa\to +\infty}\sup_{|t|\le T}
   \,\kappa\Lambda(z^\kappa(t), \zeta^\kappa(t))$ exists for all $T>0$, 
and it is given by
\begin{eqnarray*} 
   Q(T) & := & \lim_{\kappa\to +\infty}\sup_{|t|\le T}
   \,\kappa\Lambda(z^\kappa(t), \zeta^\kappa(t))
   \\ & \; = & \lim_{\kappa\to +\infty}\sup_{|t|\le T}
\Big[H_0(z^\kappa(0))-H_0(z^\kappa(t))\Big]
   =\sup_{|t|\le T}\Big[H_0(z(0))-H_0(z(t))\Big]\,.
\end{eqnarray*} 
On the other hand the equation of motion and \eqref{coer2} imply that
\begin{eqnarray*} 
   \Bigl|\,   H_0(z^\kappa(t))-H_0(z^\kappa(0))\,\Bigr|
   & = & \Bigl|\int_0^t\frac{d}{ds}\,[H_0(z^\kappa(s))]\,ds\Bigr|
   =\Bigl|\int_0^t \langle DH_0(z^\kappa(s)), \dot{z}^\kappa(s)\rangle\,ds\,\Bigr|
   \\ & = & \Bigl|\,\kappa\,\int_0^t\langle DH_0(z^\kappa(s)),
    \Pi_1\,X_{\Lambda}(z^\kappa(s), \zeta^\kappa(s))\rangle\,ds\,\Bigr| \\
    & \leq &
    C\int_0^t\sup_{|s'|\le s}\kappa\Lambda(z^\kappa(s'), \zeta^\kappa(s'))\,ds\,,
\end{eqnarray*} 
from which we deduce the Gronwall inequality $Q(T)\le C\int_0^T Q(s)\,ds$ 
in the limit $\kappa\to +\infty$. Therefore  $Q(T)=0$ for all $T\geq 0$,
and hence \eqref{urth} holds. It then follows from \eqref{coer2},
\eqref{intf} and assumption (i) that $z(t)=z(0)+\int_0^t\,\Pi_1\, X_{H_0}(z(s),0)\,ds$, 
i.e.~the curve $t\mapsto z(t)$ is the integral curve of the Hamiltonian
$H(z, 0)=H_0(z)$ starting at $z(0)$, which is unique since ${H_0}$ defines a 
Lipschitz continuous Hamiltonian vector field by assumption. 
It follows from the uniqueness of this limit curve that all subsequences have a subsequence 
which converges to the same limit, and hence that 
$(z^\kappa(t),\zeta^\kappa(t))$ converges to $(z(t),0)$ without recourse to
subsequences, as asserted in the lemma. 
\qed\bigskip

\begin{remarks}
(a) The conclusion \eqref{urth0}  says in words that in the limit 
the curve is constrained to lie on the $\zeta=0$ subspace, while
\eqref{urth}  says in words that in the limit all
the energy is in the $z$ variable, and this variable evolves in a way
that conserves $H_0(z)$ - this evolution is in fact the Hamiltonian evolution
determined by $H_0(z)$.
\smallskip

\noindent
(b) Clearly the conditions on $H,{H_0},\Lambda$ only need to hold on some
open set containing the region defined in \eqref{zzet-bd}. Also,
in \eqref{coer1} the function $c_0|\,\cdot\,|^p$ could be
replaced by any function tending to $+\infty$ at $\infty$.
\end{remarks}
\bigskip

%%%%%%%%%%%%%%%%%%%%%%%%%%%%%%%%%%%%%%%%%%%%%%%%%%%%%%%%%%%%%%%%%%

\noindent
We assume that the initial values $z(0)=(x(0),y(0))\in\R^{2n}$ 
are close to the equilibrium point
$(0, 0)$ for the real analytic Hamiltonian. 
In that case we have the following theorem:

\begin{theorem}\label{gszk} 
Let $H$ be a real analytic function of the form 
\begin{equation}\label{ham2}
   H(z,\zeta)=H_0(z)+\kappa\Lambda(z, \zeta),
   \quad\mbox{with}\quad H_0(z)=\langle\alpha, I(z)\rangle+\frac{1}{2}\,\langle AI(z), I(z)\rangle+f(z),
\end{equation}
such that $\langle AI, I\rangle\geq\frac{1}{M}\,|I|^2$ 
and $f$ is real analytic so that $f(z)={\cal O}(z^5)$ as $|z|\to 0$, and also
verifying \eqref{coer0}-\eqref{coer2}.
Fix $a\in ]0, \frac{1}{1+3n}[$. 
Then there exist positive numbers $K, k$ (depending on $a$, $n$, $\alpha$, $M$ 
and $\|A\|$) and $\theta_0$ (depending on $a$, $n$, $\alpha$, $M$, $\|A\|$ and $f$) 
with the following properties. If $t\mapsto (z(t), \zeta(t))$ is an integral curve of $H$
and $\I(0)=I(z(0))$ is such that $|\I(0)|=\theta^2$ 
for some $0<\theta\le\theta_0$, then $\I(t)=I(z(t))$ satisfies
\begin{equation}\label{gszs}
   |\I(t)-\I(0)|\leq K\theta^{2+a}
   \quad\hbox{for}\quad |t|\le e^{\frac{k}{\theta^a}}\,,
\end{equation}
for sufficiently large $\kappa$ (depending on the initial conditions
and $\theta$). 
\end{theorem}

\begin{remark} Strictly speaking to match the notation in 
lemma \ref{constr} the integral curve for $H$ 
should be written as $t\mapsto (z^\kappa(t), \zeta^\kappa(t))$, 
but we drop the additional superscript to simplify the notation. 
\end{remark}

\noindent
{\bf Beginning of proof of theorem \ref{gszk}\,} 
Following \cite{lochak1} this will be deduced from three facts:
\begin{itemize}
\item[{\rm (a)}] periodic orbits are dense in a neighbourhood
of the fixed point, 
\item[{\rm (b)}] motion in a neighbourhood of a periodic orbit
satisfies long-time stability estimates, on account of the normal
form lemma \ref{norm}, and
\item[{\rm (c)}] a priori control of the effect of the transverse component
$\zeta$ is provided by \eqref{urth}.
\end{itemize}
To begin with, since $z(0)\in\R^{2n}$ has real components, 
note that $|z(0)|^2=\sum_{j=1}^n |z_j(0)|^2=2\sum_{j=1}^n |\I_j(0)|=2|\I(0)|=2\theta^2$. 
In what follows the parameter $\theta$ will be used as a book-keeping device, 
i.e.~all quantities which need to be controlled 
will be controlled in terms of $\theta$. We are going to apply the normalization
lemma \ref{norm} to $H_0=H_0(z)$, i.e.~averaging will be performed in the $z$ variable only. 
Therefore we make the following modification of the notation defined in the introduction:

\smallskip
{\it Throughout this proof only we write 
${\cal D}_{\underline{r}}={\cal D}_{\,r_1,\,r_2}$ and
${\cal D}_{r_1,\,r_2}=\{z\in\C^{2n}: |I-{I^0}|<r_1, |z|<r_2\}$
and drop the third component from the definition of the corresponding norms
${|\,\cdot\,|}_{\underline{r}}$.}
\smallskip

First we apply corollary \ref{pds} with $I$ replaced by $\I(0)$ 
and $g=0$ in (\ref{H0def}) below. Then $\Omega(I)=\alpha+AI$ and there exist 
$K_1>0$ (depending on $\alpha$ and $A$) and $\theta_1>0$
(depending on $\alpha$, $A$, $a$ and $n$) such that the following holds.
If $|\I(0)|=\theta^2$ for some $0<\theta\le\theta_1$,
then there are ${I^0}\in\R^n$ and $\tau>0$ satisfying
\begin{enumerate}
\item
${|\I(0)-{I^0}|}_\infty\leq K_1\,\frac{\theta^{2+a}}{\tau}$, and
\item
$\pi\le\tau\le 4\pi\theta^{-a(n-1)}$
\end{enumerate}
and such that $\omega^0=\alpha+A{I^0}$ is $\tau/\theta^2$-periodic, 
i.e.~$T\omega^0\in 2\pi\Z^n$ for $T=\tau/\theta^2$. We will call this 
orbit the {\em approximating periodic orbit}.
Up to a constant, which does not affect the flow, we rewrite $H_0$ from (\ref{ham2}) as 
\begin{equation}\label{hamm}
   H_0(z)=\langle\omega^0, I\rangle+\frac{1}{2}\,
   \langle A(I-{I^0}), I-{I^0}\rangle+f(z).
\end{equation}
We will now apply the following result on stability in a neighbourhood
of periodic orbits:

\begin{lemma}[Local Stability]\label{lsz}
Consider the Hamiltonian $H$ from (\ref{ham2}). 
Assume also that $H_0$ is written as in \eqref{hamm},
with $f$ real analytic on an open neighbourhood 
of $\overline{{\cal D}}_{\,3\underline{r}}$
with ${|f|}_{\,3\underline{r}}\le\eps$ and $r_1,r_2>0$ such that
\beq
   r_1<\frac{1}{4}\,r_2^2,
   \quad\eps M<\frac{r_1^2}{2200},\quad \hbox{and}\quad |{I^0}|<\frac{r_2^2}{16}. 
   \label{ass1}
\eeq
Assume further that $\omega^0\in\R^n$ is such that $T\omega^0\in 2\pi\Z^n$, 
and that for some $m\in\N$ and $l_2>0$
\beq
54 m\|A\|r_1T\le\frac{1}{4},
\quad m^2 \eps T<\frac{l_2 r_1^2}{r_2^2}\,.
   \label{ass2}
\eeq
Let $t\mapsto (z(t),\zeta(t))$ be an integral curve for $H$ whose initial
data $(z(0),\zeta(0))$  are such that 
\begin{equation}\label{urthhere}
   \kappa\,\Lambda(z(0),\zeta(0))\le\frac{r_1^2}{360 M}\,,\quad
\hbox{
and}
\quad |\I(0)-{I^0}|\le l_1 r_1\,.
\end{equation} 
Then, for $l_1,l_2$ sufficiently small, there holds
\begin{equation}\label{conc}
   |\I(t)-{I^0}|<r_1\quad\hbox{for}\quad |t|\le t_\ast, 
\end{equation}
where $t_\ast>0$ is any time such that 
\begin{equation}\label{anyti} 
t_\ast\le\frac{3\cdot 2^m r_1}{50|\omega^0|r_2^2}
   \quad\mbox{and}\quad t_\ast\kappa\max_{|t|\leq t_\ast}\Lambda(z(t),\zeta(t))
   \leq\frac{8\eps}{5r_2 c_2(r_2)|\omega^0|}\,,
\end{equation} 
with $c_2$ from (\ref{coer2}). To be precise, the following choices 
for $l_1,l_2$ will suffice:
\beq
\quad l_1=\min\Big\{\frac{1}{4}, \frac{1}{5\sqrt{M\|A\|}}\Big\},
   \quad l_2=\min\Big\{\frac{1}{2592}, \frac{1}{120\sqrt{M\|A\|}}\Big\}\,.
   \label{ass3}
\eeq
\end{lemma}
{\bf Proof of lemma \ref{lsz}\,} As already stated we apply
the normal form lemma \ref{norm}, specialized to the
case that there is no $\zeta$ dependence, to $H_0=H_0(z)$ 
so that all the conditions
involving $r_3$ or $\kappa$ are to be disregarded, 
and also $g=0$ and $\delta=0$.
The conditions in (iv)--(vi) of that lemma are then easily seen
to be satisfied as a consequence of \eqref{ass1}, \eqref{ass2} and \eqref{ass3}.
Hence there exists a real analytic symplectic transformation
$\Psi: {\cal D}_{\,2\underline{r}}\to {\cal D}_{\,3\underline{r}}$
such that, on ${\cal D}_{{2\underline r}}$,
\[ \tilde{H}_0(\tilde{z}):=H_0\circ\Psi(\tilde{z})
=\langle\omega^0, I(\tilde{z})\rangle
   +\frac{1}{2}\,\langle A(I(\tilde{z})-{I^0}), I(\tilde{z})-{I^0}\rangle
   +\hat{g}(\tilde{z})+\hat{f}(\tilde{z}) \]
and with the properties:
\begin{itemize}
\item[{\rm (a)}] \qquad $\displaystyle {|\Psi-{\rm id}|}_{\,2\underline{r}}
\le\frac{18\,mr_2}{r_1}\,\eps T$,
\item[{\rm (b)}] \qquad ${|{\hat{g}}|}_{\,2\underline{r}}
\le 2\eps$ and $\{\hat{g}, h\}=0$,
\item[{\rm (c)}] \qquad ${|{\hat f}|}_{\,2\underline{r}}\le 2^{-m}\eps $.
\end{itemize}
The total Hamiltonian is now $\tilde{H}(\tilde{z}, \zeta)
=\tilde{H}_0(\tilde{z})+\kappa\tilde\Lambda(\tilde{z},\zeta)$ where $\Lambda(\Psi(\tilde{z}),
\zeta)=\tilde{\Lambda}(\tilde{z},\zeta)$ defines $\tilde{\Lambda}$.
Note that in order to distinguish integral curves of 
the normal form Hamiltonian $\tilde{H}$
from those of the original Hamiltonian $H$, its $z$-variables are marked by a tilde;
the relation is $z=\Psi(\tilde z)$.
We first obtain bounds for $\tilde{\I}(t)-{I^0}=I(\tilde{z}(t))-I^0$ 
for the flow of $\tilde{H}$. We will then
show that these imply (\ref{conc}) for ${\I}(t)-{I^0}=I({z}(t))-I^0$ 
with  $z(t)=\Psi({\tilde z}(t))$, using lemma \ref{pobiv} below 
to ensure that $z(t)\in  {\cal D}_{\,\underline{r}}$ 
can indeed be written thus. But for the moment we assume this and consider 
an integral curve $t\mapsto (\tilde{z}(t), \zeta(t))$ of $X_{\tilde{H}}$ 
such that $t\mapsto\tilde{z}(t)\in {\cal D}_{\,5\underline{r}/3}$. Since in general
$\{G(\tilde{I}), F(\tilde{I})\}=0$, using (b) we obtain for $\tilde{h}(t)
=h(\tilde{z}(t))=\langle\omega^0, \tilde{I}(\tilde{z}(t))\rangle$ the relation
\[ \frac{d\tilde{h}}{dt}=\langle D{h}, X_{\tilde{H}}\rangle=
\{{h}, \tilde{H}\}=\{{h}, \hat{f}+\kappa\tilde{\Lambda}\}
=\langle D{h}, X_{\hat f}\rangle
+\kappa\langle D{h}, \Pi_1 X_{\tilde\Lambda}\rangle\,. 
\]
Next observe that
\begin{equation}\label{Cauprep}
   \tilde{z}\in {\cal D}_{{5\underline r/3}}\quad\mbox{and}\quad
   |\tilde{w}-\tilde{z}|\le\frac{r_1}{10\,r_2}\quad\Longrightarrow\quad\tilde{w}\in {\cal D}_{{2\underline r}},
\end{equation}
since $|\tilde{w}|\le |\tilde{w}-\tilde{z}|+|\tilde{z}|<r_1/10r_2+5r_2/3<2r_2$ 
by the condition $r_1<r_2^2/4$ in (\ref{ass1});
using in addition (\ref{Iwz}) we obtain 
\[ |{I}(\tilde{w})-{I^0}|\le |{I}(\tilde{w})-{I}(\tilde{z})|
   +|{I}(\tilde{z})-{I^0}|\le\frac{1}{2}\,(|\tilde{w}-\tilde{z}|+2|\tilde{z}|)|\tilde{w}-\tilde{z}|
   +\frac{5r_1}{3}<2r_1. \]  
Thus we can bound by means of Cauchy's estimate as in (\ref{cauch}):
\begin{equation}\label{bmoc}
   {|X_{\hat{f}}|}_{\,5\underline{r}/3}
   \leq\frac{10\,r_2 {|{\hat{f}}|}_{\,2\underline{r}}}{r_1}
   \leq\frac{10\,r_2 2^{-m}\eps}{r_1}\,.
\end{equation}
Also by Cauchy's estimate, (\ref{Cauprep}) and (a) we obtain
\begin{equation}\label{derbd} 
   {|D\Psi-1|}_{\,5\underline{r}/3}
   \le\frac{10\,r_2}{r_1}\,{|\Psi-{\rm id}|}_{\,2\underline{r}}
   \le\frac{180\,m r_2^2}{r_1^2}\,\eps T
   \le\frac{180\,m^2 r_2^2}{r_1^2}\,\eps T<180\,l_2<\frac{1}{2}\,, 
\end{equation}
from which we derive the pointwise estimate
$|\Pi_1 X_{\tilde\Lambda}|\le |\frac{\partial\Lambda}{\partial z}||D\Psi|\leq 3c_2(r_2)\Lambda/2$, 
using also \eqref{coer2}. It follows that, 
for as long as $\tilde{z}(t)$ remains in ${\cal D}_{\,5\underline{r}/3}$
and $z(t)=\Psi(\tilde{z}(t))\in {\cal D}_{\,\underline{r}}$,
\beq\la{tm}
   \bigg|\frac{d\tilde{h}}{dt}(t)\bigg|
   =|\langle D\tilde{h}, X_{\hat f} +\kappa\Pi_1 X_{\tilde\Lambda}
   \rangle|\leq\frac{5}{3}\,|\omega^0|\,r_2
   \,\Bigl(\frac{10\,r_2 2^{-m}\eps}{r_1}+\frac{3\kappa}{2}\,c_2(r_2)
   \,|{\Lambda}(z(t), \zeta(t))|\Bigr)\,.
\eeq
>From the definition of $t_\ast$ we deduce that, for $|t|\le t_\ast$ as 
in \eqref{conc}-\eqref{anyti},
\[ |\tilde{h}(t)-\tilde{h}(0)|\leq\eps+\frac{5\kappa}{2}\,r_2 c_2(r_2)|\omega^0|
   \,|t|\max_{|t'|\leq |t|}\Lambda(z(t'),\zeta(t'))\leq 5\eps\,. \] 
Energy conservation 
$\tilde{H}(\tilde{z}(t), \zeta(t))=\tilde{H}(\tilde{z}(0), \zeta(0))$ together
with \eqref{coer0} 
and the convexity assumption (strict positivity of the matrix $A$) then give:
\begin{eqnarray*}
   \frac{1}{2M}\,{|\tilde{\I}(t)-{I^0}|}^2
   & \leq & \frac{1}{2}\,\|A\|\,{|\tilde{\I}(0)-{I^0}|}^2
   +|\tilde{h}(t)-\tilde{h}(0)|+2{|{\hat{g}}|}_{\,5\underline{r}/3}
   +2{|{\hat{f}}|}_{\,5\underline{r}/3}+\kappa\,\Lambda(z(0), \zeta(0))
   \\
   & \leq & \frac{1}{2}\,\|A\|\,{|\tilde{\I}(0)-{I^0}|}^2
   +5\eps+6\eps+\kappa\Lambda(z(0),\zeta(0)), 
\end{eqnarray*}
so that due to (\ref{ass1})-(\ref{ass3}) and for $|t|\le t_\ast$:  
\begin{eqnarray}\label{transI}
   {|\tilde{\I}(t)-{I^0}|}^2 
   & \le & M\|A\|\,{|\tilde{\I}(0)-{I^0}|}^2+22\,\eps M
   +2M\kappa\Lambda(z(0),\zeta(0))
   \nonumber
   \\ & \le & M\|A\|\,{|\tilde{\I}(0)-{I^0}|}^2+\frac{r_1^2}{100}
   +2M\kappa\Lambda(z(0),\zeta(0))
\end{eqnarray}
for as long as $\tilde{z}(t)\in {\cal D}_{\,5\underline{r}/3}$
and $z(t)\in {\cal D}_{\,\underline{r}}$. 

Now to deduce (\ref{conc}) it is necessary both to show that
\eqref{transI} implies the inequality in \eqref{conc}, and also to justify
the assumption that $\tilde{z}(t)\in {\cal D}_{\,5\underline{r}/3}$
and ${z}(t)\in {\cal D}_{\,\underline{r}}$ made above in deriving \eqref{transI}.
To this end suppose that $|\I(0)-{I^0}|\le l_1 r_1$
for an integral curve $t\mapsto (z(t), \zeta(t))$ of the original Hamiltonian
vector field $X_H$. Since we are considering real-valued solutions 
of the Hamiltonian equations, 
\begin{eqnarray*} 
   |z(0)|^2 & = & 2|\I(0)|\le 2(|\I(0)-{I^0}|+|{I^0}|)
   \\ & \le & 2 l_1 r_1+\frac{1}{8}\,r_2^2
   <\Big(\frac{l_1}{2}+\frac{1}{8}\Big)r_2^2\le\frac{r_2^2}{4}. 
\end{eqnarray*} 
Therefore we have $z(0)\in {\cal D}_{\,\underline{r}/2}$. 
Denote by $t_0>0$ the longest time such that $z(t)\in {\cal D}_{\,\underline{r}}$
for all $|t|\le t_0$. 

The point of the following lemma \ref{pobiv} is to show that
a sufficiently large neighbourhood of the approximating 
periodic orbit is covered by the transformation $\Psi$
(as a consequence of (a) and the various assumptions on the parameters used).
This ensures that stability information just derived for integral curves of 
the transformed Hamiltonian $\tilde H$ will imply stability information
for the integral curves of $H$ on a sufficiently large neighbourhood of 
this periodic orbit.
Here we write 
${\cal D}_{\,\underline{r}}^{({\rm real})}={\cal D}_{r_1,\,r_2}^{({\rm real})}$ where
\[ {\cal D}_{a,\,b}^{({\rm real})}=\{z\in\R^{2n}: |I(z)-{I^0}|<a, |z|<b\}, \]
and similarly we denote 
\[ B_r^{({\rm real})}(w)=\{z\in\R^{2n}: |z-w|<r\} \] 
for $r>0$ and $w\in\R^{2n}$. 

\begin{lemma}\label{pobiv}
Under the hypotheses of lemma \ref{lsz}, $\Psi$ satisfies
$\Psi({\cal D}_{\,5\underline{r}/3}^{({\rm real})})
\supset {\cal D}_{\,\underline{r}}^{({\rm real})}$.
\end{lemma}
{\bf Proof of lemma \ref{pobiv}\,} According to (a) and (\ref{derbd}) 
we have ${|\Psi-{\rm id}|}_{\,2\underline{r}}\le\frac{18\,mr_2}{r_1}\,\eps T=:\mu$ 
and ${|D\Psi-1|}_{\,5\underline{r}/3}<1/2$. 
Hence $D\Psi(z)$ is invertible for every $z\in {\cal D}_{\,5\underline{r}/3}$, and accordingly 
$\Psi: {\cal D}_{\,5\underline{r}/3}\to \Psi({\cal D}_{\,5\underline{r}/3})=:{\cal W}$ is a real-analytic 
diffeomorphism such that $\|D\Psi^{-1}(w)\|\le 2$ for $w\in {\cal W}$.
Now fix $w\in {\cal D}_{\,\underline{r}}^{({\rm real})}$. 
Then $B_\delta^{({\rm real})}(w)\subset {\cal D}_{\,3\underline{r}/2}^{({\rm real})}$ 
for $\delta=\frac{r_1}{4\,r_2}$, as can be shown 
using  $r_1<r_2^2/4$ and (\ref{Iwz}),
analogously to (\ref{Cauprep}).
Furthermore, for $w\in {\cal D}_{\,\underline{r}}^{({\rm real})}$, 
\[ |w-\Psi(w)|\le\mu<\frac{\delta}{2} \] 
due to $\frac{18\,mr_2}{r_1}\,\eps T\le\frac{4\cdot 18\,m^2 r_2^2}{r_1^2}
\,\eps T\times\frac{r_1}{4r_2}\le 72\,l_2\delta<\frac{\delta}{2}$. 
In other words, we have $w\in B_{\delta/2}^{({\rm real})}(\Psi(w))$. 
Next we apply lemma \ref{LipOM} below and use the fact that $\Psi$ is 
real on real vectors,
to deduce that  
\[ \Psi\Big(\,\overline{B_\delta^{({\rm real})}(w)}\,\Big)
   \supset\overline{B_{\delta/2}^{({\rm real})}(\Psi(w))}. \] 
To summarize, for fixed $w\in {\cal D}_{\,\underline{r}}^{({\rm real})}$ we obtain 
\[ w\in B_{\delta/2}^{({\rm real})}(\Psi(w))
   \subset\Psi\Big(\,\overline{B_\delta^{({\rm real})}(w)}\,\Big)
   \subset\Psi\Big(\,\overline{{\cal D}_{\,3\underline{r}/2}^{({\rm real})}}\,\Big)
   \subset\Psi({\cal D}_{\,5\underline{r}/3}^{({\rm real})}), \]  
and this concludes the proof of lemma \ref{pobiv}. 
\qed\bigskip

\noindent
{\bf Continuation of the proof of lemma \ref{lsz}}
Due to Lemma \ref{pobiv} we may write 
$z(t)=\Psi(\tilde{z}(t))$ for $|t|\le t_0$
with an integral curve $t\mapsto (\tilde{z}(t), \zeta(t))$ 
of $X_{\tilde{H}}$ such that $t\mapsto\tilde{z}(t)
\in {\cal D}_{\,5\underline{r}/3}^{({\rm real})}$. 
Then by (a) and (\ref{ass1})-(\ref{ass3}),
\begin{eqnarray}\label{diff0}
   |\tilde{\I}(0)-{I^0}| & \le & |\tilde{\I}(0)-\I(0)|+|\I(0)-{I^0}|
   \le\frac{1}{2}\,\Big(|\tilde{z}(0)|+|z(0)|\Big)\,|\tilde{z}(0)-\Psi(\tilde{z}(0))|+l_1 r_1
   \nonumber
   \\ & \le & \frac{1}{2}\,\Big(\frac{5r_2}{3}+r_2\Big)\,\frac{18\,mr_2}{r_1}
\,\eps T+l_1 r_1
   \le\frac{24\,m^2 r_2^2}{r_1}\,\eps T+l_1 r_1\le (24\,l_2+l_1)\,r_1.
\end{eqnarray}
Then we can apply (\ref{transI}) and use $t_\ast\le T$ together with (\ref{urthhere}) to obtain 
\begin{eqnarray*} 
   {|\tilde{\I}(t)-{I^0}|^2}
   & \le & {M\|A\|}\,|\tilde{\I}(0)-{I^0}|^2+\frac{r_1^2}{100}
   +2M\kappa\,\Lambda(z(0),\zeta(0))
   \\ & \le & \Big({M\|A\|}\,(24\,l_2+l_1)^2+\frac{1}{100}\Big)\,r_1^2
   +\frac{r_1^2}{180}   
   \\ & \le & \Big(\frac{4}{100}+\frac{1}{100}\Big)\,r_1^2+\frac{r_1^2}{180}
   \\ & = & \frac{r_1^2}{20}+\frac{r_1^2}{180}<\Big(\frac{r_1}{4}\Big)^2
\end{eqnarray*} 
for $|t|\le\min\{t_\ast, t_0\}$, since $l_2$, $l_1$ are such that
${M\|A\|}\,(24\,l_2+l_1)^2\leq\frac{4}{100}=(\frac{1}{5})^2$.
In the same manner as for (\ref{diff0}) this in turn leads to
\begin{eqnarray}\label{Idiffs}
   |\I(t)-{I^0}| & \le & |\I(t)-\tilde{\I}(t)|+|\tilde{\I}(t)-{I^0}|
   \nonumber
   \\ & \le & \frac{1}{2}\,\Big(|\tilde{z}(t)|+|z(t)|\Big)
   \,|\tilde{z}(t)-\Psi(\tilde{z}(t))|+\frac{r_1}{4}
   \nonumber
   \\ & \le & \Big(24\,l_2+\frac{1}{4}\Big) r_1<\frac{r_1}{2}
\end{eqnarray}
for $|t|\le\min\{t_\ast, t_0\}$, due to $24l_2<\frac{1}{4}$. 
Since also $r_1<r_2^2/4$, this implies that for such times
\[ |z(t)|^2 =  2|\I(t)|\le 2\Big(|\I(t)-{I^0}|+|{I^0}|\Big)
   \le 2\Big[\frac{r_2^2}{8}
   +\frac{r_2^2}{16}\Big]\,<r_2^2\,.
\]
Hence we see that $\min\{t_\ast, t_0\}<t_0$, 
or in other words $\min\{t_\ast, t_0\}=t_\ast$.
Thus (\ref{conc}) is a consequence of (\ref{Idiffs}). \qed\bigskip

\noindent
{\bf Completion of proof of theorem \ref{gszk}\,}
We now aim to show that the stability bound (\ref{conc}), applied
in the neighbourhood of the approximating periodic orbit obtained
prior to lemma \ref{lsz}, implies (\ref{gszs}).
Since $f$ vanishes to fifth order we take $r_2=8\theta$ and $\eps=C_1\theta^5$ 
to ensure that ${|f|}_{\,3\underline{r}}\le\sup\,\{|f(z)|: |z|\le 3r_2\}
\le C_0(3r_2)^5\le\eps$, where $C_1={24}^5\,C_0$ has to be chosen large enough (depending on $f$). 
In addition, let 
\begin{equation}
\label{impdef}
m=\delta\,[\theta^{-a}]\quad\mbox{and}\quad r_1=\frac{L\theta^{2+a}}{\tau}, 
\end{equation}
where $\delta, L>0$ will be fixed below; recall that the period
of the approximating periodic orbit is $T=\tau/\theta^2$. We will now verify
that having fixed $l_1$, $l_2$ satisfying (\ref{ass3}), 
the conditions (\ref{ass1})-(\ref{ass2}) can be made to hold
by making $\theta$ sufficiently small and choosing $\delta$, $L$ appropriately. 
To start with
\[ \frac{r_1}{r_2^2}=\frac{L\theta^a}{64\tau}\le\frac{L\theta^a}{64\pi} \] 
by (ii), and hence the first condition of (\ref{ass1}) holds if $\theta$ is small enough. 
In addition, 
\[ m r_1 T=\delta\,[\theta^{-a}]\,\frac{L\theta^{2+a}}{\tau}\,\frac{\tau}{\theta^2}\le\delta L \] 
whence we need to have 
\begin{equation}\label{delL1} 
   \delta L\le\frac{1}{216\,\|A\|}
\end{equation} 
to validate the first condition of (\ref{ass2}). Next, 
\[ \frac{\eps}{r_1^2}=\frac{C_1\theta^5\tau^2}{L^2\theta^{4+2a}}
   \le\frac{16\pi^2 C_1}{L^2}\,\theta^{1-2an} \] 
by (ii) shows that we can fulfil the second condition of (\ref{ass1}) 
for $\theta$ sufficiently small, due to $a<\frac{1}{2n}$. 
Concerning the condition on $|{I^0}|$ in (\ref{ass1}), here 
\[ |{I^0}|\le |\I(0)-{I^0}|+|\I(0)|\le nK_1\,\frac{\theta^{2+a}}{\tau}+\theta^2
   \le\Big(nK_1\,\frac{\theta^a}{\pi}+1\Big)\theta^2\le 2\theta^2 \]
by (i) and (ii) for $\theta$ small enough. Hence 
\[ \frac{|{I^0}|}{r_2^2}\le \frac{2\theta^2}{64\theta^2}<\frac{1}{16}, \] 
and thus all of (\ref{ass1}) is verified, provided that (\ref{delL1}) can be ensured. 
To establish the second condition of (\ref{urthhere}), note that 
\[ \frac{|\I(0)-{I^0}|}{r_1}\le nK_1\,\frac{\theta^{2+a}}{\tau}\,\frac{\tau}{L\theta^{2+a}}
   =\frac{nK_1}{L} \] 
by (i). Accordingly, we need to have 
\begin{equation}\label{delL2} 
   \frac{nK_1}{L}\le l_1
\end{equation} 
for $l_1$ from (\ref{ass3}). For the last condition of (\ref{ass2}) finally 
\[ \frac{m^2\eps T r_2^2}{r_1^2}
   \le\frac{\delta^2\theta^{-2a}C_1\theta^5\tau^3\cdot 64\theta^2}{\theta^2 L^2\theta^{4+2a}}
   =\frac{64\,C_1\delta^2}{L^2}\,\theta^{1-4a}\,\tau^3
   \le\frac{4096\,\pi^3C_1\delta^2}{L^2}\,\theta^{1-a(1+3n)} \]  
by (ii). Since $a<\frac{1}{1+3n}$, the right-hand side is smaller than $l_2$ from (\ref{ass3}), 
if $\theta$ is sufficiently small. Altogether, (\ref{ass1}) and (\ref{ass2}) 
will be satisfied, provided that (\ref{delL1}) and (\ref{delL2}) hold. 
This can be achieved by explicitly taking 
\[ L=\frac{nK_1}{l_1}\quad\mbox{and}\quad\delta=\frac{1}{216\,\|A\| L}. \] 
We thus have shown so far that there is $\theta_0>0$ 
(depending on the quantities as stated in the theorem) such that for 
$0<\theta\le\theta_0$ 
the assumptions (\ref{ass1}) and (\ref{ass2}) from Lemma \ref{lsz} 
hold, as does the second condition of (\ref{urthhere}).
Now fix $0<\theta\le\theta_0$ and put 
$t_\ast=\frac{3\cdot 2^m r_1}{50|\omega^0|r_2^2}$ 
(depending on $\theta$). Then (\ref{urth}) from Lemma \ref{constr} ensures that 
\[ \lim_{\kappa\to +\infty}\max_{|t|\leq t_\ast}\,\kappa\,
\Lambda (z(t),\zeta(t))=0\,; \] 
(recall that $(z(t),\zeta(t))=(z^\kappa(t),\zeta^\kappa(t))$ 
in the notation of lemma \ref{constr}). In particular, 
the first condition of (\ref{urthhere}) and the third condition 
from (\ref{anyti}) will be satisfied, 
if $\kappa\ge\kappa_0$ for an appropriate $\kappa_0=\kappa_0(\theta)>0$ 
depending on the initial data and $\theta$. 
Therefore Lemma \ref{lsz} applies and we deduce from (\ref{conc}) that 
$|\I(t)-{I^0}|<r_1$ for $|t|\le t_\ast$. Now combine this with (i) to bound,
for $|t|\le t_\ast$,
$$
|\I(t)-\I(0)|\leq |\I(t)-{I^0}|+|{I^0}-\I(0)|\leq r_1+n K_1\,\frac{\theta^{2+a}}{\tau}
\le\frac{2L}{\pi}\,\theta^{2+a}=K\theta^{2+a}\,,
$$
where we have defined 
$K=\frac{2L}{\pi}$. (The penultimate inequality holds since 
$r_1=\frac{L\theta^{2+a}}{\tau}$ and $L\ge nK_1$). 

It remains to observe that by (ii), 
\[ t_\ast=\frac{3\cdot 2^m L\theta^a}{3200\,|\omega^0|\tau}
   \ge\frac{3\cdot 2^m L\theta^{an}}{12800\,\pi |\omega^0|}, \] 
so that with $B=\frac{3L}{12800\,\pi |\omega^0|}$ and $m$ as in
(\ref{impdef}) and for $\theta$ small enough 
(reducing $\theta_0$ further if necessary) 
\[ \ln t_\ast\ge\ln B+an\ln\theta+(\ln 2)\delta\,[\theta^{-a}]
   \ge k\theta^{-a} \] 
for any $k<({\ln 2}){\delta}$, and in particular for 
$k=\frac{\ln 2}{2}\,\delta$, completing the proof of \eqref{gszs} and the 
theorem.
\qed\bigskip

As already remarked theorem \ref{gszk} does not provide quantitative
information on the domains on which the bound \eqref{gszs} holds, 
only the assurance that it holds for sufficiently large $\kappa$. However
when the nonlinear interaction has a special structure it is
possible to extract precise information on the domains as we 
now explain. We assume that there are additional smooth functions 
$J_k$, $k=1,\dots l$, of $\zeta\in\R^{2N}$ which all Poisson commute with $\Lambda$:
$$
\{J_k, \Lambda\}=0\quad\mbox{for}\quad k=1,\dots, l, 
$$
and that for all $R>0$ there exists $c_3(R)>0$ such that
\begin{equation}\label{coer2p}
\left|\frac{\partial \Lambda}{\partial z}(z,\zeta)\right |
\leq c_3(R)\sum_{k=1}^l |J_k(\zeta)|\quad\hbox{for}\quad |z|\leq R\,.
\end{equation}

Then we have the following quantitative version of theorem \ref{gszk}:

\begin{theorem}\label{quantivers} 
Let $H$ be a real analytic function of the form \eqref{ham2}
such that $\langle AI, I\rangle\geq\frac{1}{M}\,|I|^2$ 
and $f$ is real analytic so that $f(z)={\cal O}(z^5)$ as $|z|\to 0$, and also
verifying \eqref{coer0}-\eqref{coer1} and \eqref{coer2p}.
Fix $a\in ]0, \frac{1}{1+3n}[$. 
Then there exist positive numbers $K, k$ (depending on $a$, $n$, $\alpha$, $M$ 
and $\|A\|$) and $\theta_0$ (depending on $a$, $n$, $\alpha$, $M$, $\|A\|$ and $f$) 
with the following properties. If $t\mapsto (z(t), \zeta(t))$ is an integral curve of $H$
and $\I(0)=I(z(0))$ is such that $|\I(0)|=\theta^2$ 
for some $0<\theta\le\theta_0$, then $\I(t)=I(z(t))$ satisfies
\[ |\I(t)-\I(0)|\leq K\theta^{2+a}
   \quad\hbox{for}\quad |t|\le e^{\frac{k}{\theta^a}} \] 
for initial data such that
\beq\la{holds}
   \sum_{k=1}^l |J_k(\zeta(0))|\leq
   \frac{\theta^4 e^{-\frac{k}{\theta^a}}}{\kappa}\quad\hbox{and}\quad
\kappa\Lambda(z(0),\zeta(0))\leq\frac{L^2\theta^{4+2an}}{(4\pi)^2360 M}\,.
\eeq
\end{theorem}
\noindent
{\bf Proof\,}
The proof is almost entirely the same as the proof of theorem
\ref{gszk} except for two points: 

firstly, the condition
\eqref{urthhere} required to apply lemma \ref{lsz} is an explicit
consequence of the second inequality in \eqref{holds}; and

secondly, to bound $\frac{d\tilde{h}}{dt}$
the estimate \eqref{tm} is now replaced by
\[ \bigg|\frac{d\tilde{h}}{dt}(t)\bigg|
   =|\langle D\tilde{h}, X_{\hat f} +\kappa\Pi_1 X_{\tilde\Lambda}
   \rangle|\leq\frac{5}{3}\,|\omega^0|\,r_2
   \,\Bigl(\frac{10\,r_2 2^{-m}\eps}{r_1}+\frac{3\kappa}{2}\,c_3(r_2)
   \,\sum_{k=1}^l |J_k(\zeta(0))|\Bigr)\,. \] 
(The fact that this holds with the $J_k$ evaluated at $\zeta(0)$ is a consequence
of the assumption that they Poisson commute with $\Lambda$ and so are
constants of motion.) To ensure that $|\tilde{h}(t)-\tilde{h}(0)|\leq 5\eps$ 
for $|t|\le e^{\frac{k}{\theta^a}}$ we require 
\[ e^{\frac{k}{\theta^a}}
   \frac{5\kappa}{2}\,r_2 c_3(r_2)|\omega^0|\,\sum_{l=1}^k |J_k(\zeta(0))|\leq 4\eps \] 
which, using the definitions $r_2=8\theta$ 
and $\eps=C_1\theta^5$ from the paragraph preceding \eqref{impdef}, is a consequence of \eqref{holds}
for sufficiently small $\theta$ (if necessary modifying some constants and using $c_3(8\theta)\le C$ 
for $\theta\le 1$). \qed

%%%%%%%%%%%%%%%%%%%%%%%%%%%%%%%%%%%%%%%%%%%%%%%%%%%%%%%%%%%%%%%%%%%%%%%%%%%%%%%%%%%%%%%%

\setcounter{equation}{0}

\section{Nekhoroshev stability in the case of small $\kappa$}
\label{nssk}

In this section $t\mapsto\bigl(z(t),\zeta(t)\bigr)\in\R^{2n}\times\R^{2N}$ 
is an integral curve of the real-analytic Hamiltonian 
\beq\la{ra}
H(z, \zeta)=\langle\alpha, I(z)\rangle
   +\frac{1}{2}\,\langle AI(z), I(z)\rangle
   +f_\kappa(z, \zeta)+\kappa\Lambda(\zeta). 
\eeq
We will consider the case that $f_\kappa$ is allowed to depend on $\kappa$ and satisfies 
\begin{equation}\label{fifth}
   |f_\kappa(z, \zeta)|\leq C_0\bigl(|z|^5+|\zeta|^2|z|^4+\kappa |\zeta|^2|z|\bigr)
\end{equation}
in a sufficiently large neighbourhood of the origin.
In addition we will always assume that
\begin{equation}\label{spos}
   \langle AI, I\rangle\geq\frac{1}{M}\,|I|^2\quad\hbox{and}\quad
   \Lambda(\zeta)\geq\frac{|\zeta|^2}{2}\,,
\end{equation}
and
\begin{equation}\label{spos2}
   \Lambda(\zeta)\le\frac{C_\Lambda}{2}\,|\zeta|^2
   \quad\hbox{and}\quad
   |D\Lambda(\zeta)|\le C_\Lambda |\zeta|.
\end{equation}
(These conditions are all understood
to hold on some open set in $\R^{2n}\times\R^{2N}$ 
in which the integral curve lies.)

We will prove that exponential stability estimates like
(\ref{gszs}) hold for the projected motion in the $z$-plane, together with long time
bounds for $\zeta(t)$, as long as $\kappa$ is sufficiently small.

\begin{theorem}\label{agszsk} 
Let $t\mapsto\bigl(z(t),\zeta(t)\bigr)$  be
an integral curve of the real-analytic Hamiltonian $H$ verifying
\eqref{ra}-\eqref{spos2}.
There exist constants $\kappa_0, k>0$ 
and $p_1\in ]0, 1[$, $q_1>1$, $p_2\in ]p_1, 1[$
with the following properties. If $0<\kappa\le\kappa_0$ 
and if the initial data are such that 
$\I(0)=I(z(0))$ and $\mathbf{\Lambda}(0)=\Lambda(\zeta(0))$ satisfy
\[ |\I(0)|=O(\kappa^{p_1})\,,\quad
   \kappa\mathbf{\Lambda}(0)=O(\kappa^{q_1})\,, \] 
then the quantities $\I(t)=I(z(t))$ and $\mathbf{\Lambda}(t)=\Lambda(z(t))$ satisfy
\[ |\I(t)-\I(0)|=O(\kappa^{p_2})\quad\hbox{and}\quad
   \kappa\mathbf{\Lambda}(t)=O(\kappa^{2p_2})
   \quad\hbox{for}\quad |t|\le e^{\frac{k}{\kappa^{q_2}}}\,,
\] 
where $q_2=p_2-p_1$ and $2p_2>1$. 
All of the exponents and implicit constants  
are independent of $N$.
\end{theorem}

This theorem will follow from:

\begin{theorem}\label{gszsk} 
Let $t\mapsto\bigl(z(t),\zeta(t)\bigr)$  be
an integral curve of the real-analytic Hamiltonian $H$ verifying
\eqref{ra}-\eqref{spos2}.
Fix $a\in ]0, \min\{\frac{1}{4(n-1)}, \frac{1}{1+3n}\}[$. 
Then there exist positive numbers $C_E$, $\theta_0<1$, $K$, $k$ 
with the following properties. If the initial data are such that 
$\I(0)=I(z(0))$ and $\mathbf{\Lambda}(0)=\Lambda(\zeta(0))$ satisfy
\begin{equation}\label{agsz}
   |\I(0)|\leq\theta^2\,,\quad
   \kappa\mathbf{\Lambda}(0)\leq C_E\theta^{4+2an}
   \quad\hbox{and}\quad
   \kappa = \theta^{2+2a(2n-1)}\,,
\end{equation}
for $0<\theta\le\theta_0$, 
then $\I(t)=I(z(t))$ and $\mathbf{\Lambda}(t)=\Lambda(z(t))$ satisfy
\begin{equation}\label{gszslk}
   |\I(t)-\I(0)|\leq K\theta^{2+a}\quad\hbox{and}\quad
   \kappa\mathbf{\Lambda}(t)\leq K\theta^{4+2a}
   \quad\hbox{for}\quad |t|\le e^{\frac{k}{\theta^a}}.
\end{equation}
The numbers $C_E$, $\theta_0$, $K$, $k$ depend on 
$a$, $n$, $\|A\|$, $C_0$, $M$, $C_\Lambda$, but not on $N$.
\end{theorem}

\begin{remarks} (a) Theorem \ref{agszsk} is a direct consequence 
of theorem \ref{gszsk}: it suffices to take $a\in ]0, \min\{\frac{1}{4(n-1)}, \frac{1}{1+3n}\}[$ 
and define 
\begin{eqnarray*} 
   & & p_1=\frac{2}{2+2a(2n-1)},\quad
   q_1=\frac{4+2an}{2+2a(2n-1)},
   \\ & & p_2=\frac{2+a}{2+2a(2n-1)},\quad
   q_2=\frac{a}{2+2a(2n-1)}.
\end{eqnarray*}
Notice that $a<\frac{1}{1+3n}<\frac{1}{2(n-1)}$ 
implies that
$p_2$ thus defined satisfies $2p_2>1$. 
\smallskip 

\noindent
(b) It will become apparent from the proof that the result is valid under
conditions on $f$ more general than (\ref{fifth}). The crucial thing
is that on an appropriate neighbourhood $f$ is bounded by a number $\eps$
satisfying the conditions in (\ref{ass1-2nd}) and (\ref{ass2-2nd})
and satisfying the scaling relations in the last section of the proof.
\end{remarks}

\noindent
{\bf Beginning of proof of theorem \ref{gszsk}\,} 
We follow the same basic strategy as in the proof of theorem \ref{gszk},
and start in identical fashion by introducing an approximating periodic orbit
by corollary \ref{pds} (with $I$ replaced by $\I(0)$ and $g=0$). This provides
a frequency vector
$\omega^0=\alpha+A{I^0}$ which is $\tau/\theta^2$-periodic, 
i.e.~$T\omega^0\in 2\pi\Z^n$ for $T=\tau/\theta^2$, such that:
\begin{enumerate}
\item
$\max\limits_{{{1\leq j\le n}}}{|\I_j(0)-{I_j^0}|}
\leq C_A\,\frac{\theta^{2+a}}{\tau}$, and
\item
$\pi\le\tau\le 4\pi\theta^{-a(n-1)}$. 
\end{enumerate}
This is a periodic orbit for the 
unperturbed $z$ part of the motion. The number $C_A$
is just a bound for the inverse of the map
$I\mapsto\alpha+A I$ and depends on $M$, $n$. 
Up to a constant, which does not affect the flow, we rewrite $H$ as:
\begin{equation}
\label{chamwg}
H(z, \zeta)=\langle\omega^0, I\rangle
+\frac{1}{2}\,\langle A(I-{I^0}), I-{I^0}\rangle
+f_\kappa(z, \zeta)+\kappa\Lambda(\zeta)\,. 
\end{equation}
We will now apply the following result on stability in a neighbourhood
of periodic orbits, which is the analogue of lemma \ref{lsz}:

\begin{lemma}[Stability in a neighbourhood of a periodic orbit]\label{lszsk}
Assume that $\omega^0\in\R^n$ is such that $T\omega^0\in 2\pi\Z^n$.
Consider a Hamiltonian of the form (\ref{chamwg}), verifying (\ref{spos}) 
and (\ref{spos2}), 
which is real analytic on an open neighborhood 
of $\overline{{\cal D}}_{\,3\underline{r}}$ so that ${|f_\kappa|}_{\,3\underline{r}}\le\eps$ 
and with $r_1, r_2, r_3>0$ such that
\begin{equation}
   r_1<\min\Big\{\frac{1}{4}\,r_2^2\,,\,2r_2r_3\Big\}\,,
   \quad \eps M<l_0{r_1^2}\,,\quad |{I^0}|<\frac{r_2^2}{16}\,,
   \quad \hbox{and}\quad r_1^2\le 4\kappa M r_3^2\,,
   \label{ass1-2nd}
\end{equation}
for some positive $l_0$.
Assume further $m$ is a positive integer such that
\begin{equation}
54 m\|A\|r_1T\le\frac{1}{6}\,,\quad m^2 \eps T<\frac{l_2 r_1^2}{r_2^2}\,,
\quad \hbox{and}\quad C_\Lambda\kappa m T\leq\frac{r_1}{54 r_2r_3}
\,.\label{ass2-2nd}   
\end{equation}
Then for initial data satisfying
\begin{equation}\label{ass4-2nd}
   |\I(0)-{I^0}|\le l_1 r_1,
   \quad \kappa\mathbf{\Lambda}(0)
   \le\frac{r_1^2}{200M}
\end{equation}
and with $l_0, l_1, l_2>0$ sufficiently small (depending only on $M$, $\|A\|$)
\begin{equation}\label{conc-2nd}
   |\I(t)-{I^0}|\leq r_1,
   \quad\kappa\mathbf{\Lambda}(t)\leq\frac{r_1^2}{16M},
   \quad\hbox{and}\quad
   |\zeta(t)|\leq r_3
   \quad\hbox{for}\quad |t|\le t_\ast=\frac{3\cdot 2^m r_1}{10|\omega^0|r_2^2}.
\end{equation}
To be specific the following choices for $l_0,l_1,l_2$ will suffice:
\begin{equation}
   l_0=\frac{1}{2200}\,,\quad
   l_1=\min\Big\{\frac{1}{4}, \frac{1}{20\sqrt{M\|A\|}}\Big\},
   \quad l_2=\min\Big\{\frac{1}{3888}, \frac{1}{480\sqrt{M\|A\|}}\Big\}\,.
   \label{ass3-2nd}
\end{equation}
\noindent
\end{lemma}
{\bf Proof of lemma \ref{lszsk}\,} We apply
the normal form lemma \ref{norm} with $g$ and $\delta$ set to zero:
the conditions in (iv)--(vi) of that lemma are then easily seen
to be satisfied as a consequence of (\ref{ass1-2nd})-(\ref{ass2-2nd}), 
with $l_2$ as in (\ref{ass3-2nd}).
Hence there exists a real analytic symplectic transformation
$\Psi: {\cal D}_{\,2\underline{r}}\to {\cal D}_{\,3\underline{r}}\,,$
such that on ${\cal D}_{{2\underline r}}$,
\[ \tilde{H}:=H\circ\Psi=\langle\omega^0, I(\tilde{z})\rangle
   +\frac{1}{2}\,\langle A(I(\tilde{z})-{I^0}), I(\tilde{z})-{I^0}\rangle
   +\hat{g}(\tilde{z}, \tilde{\zeta})+\hat{f}_\kappa(\tilde{z}, \tilde{\zeta})
   +\kappa\Lambda(\tilde{\zeta}) \]
and with the properties:
\begin{itemize}
\item[{\rm (a)}] \qquad $\displaystyle {|\Psi-{\rm id}|}_{\,2\underline{r}}
\le\frac{18\,mr_2}{r_1}\,\eps T=:\mu$,
\item[{\rm (b)}] \qquad ${|{\hat{g}}|}_{\,2\underline{r}}
\le 2\eps$ and $\{\hat{g}, h\}=0$ for $h=\langle\omega^0, I\rangle$, 
\item[{\rm (c)}] \qquad ${|{\hat f}_\kappa|}_{\,2\underline{r}}\le 2^{-m}\eps $.
\end{itemize}
(Notice that $\kappa$ is fixed in lemma \ref{norm}, so that lemma 
can be applied
to $f_\kappa$ depending on $\kappa$ 
and yields a new $\hat{f}_\kappa$, also depending on $\kappa$, obeying
the bound in (c)). 
The variables in the normal form Hamiltonian $\tilde{H}$
are distinguished by a tilde, and are related to the original variables 
by $(z,\zeta)=\Psi(\tilde z, \tilde{\zeta})$.
The crucial point is the small rate of change of
$\tilde{h}(t)=h(\tilde{z}(t))=\langle\omega^0, \tilde{\I}(t)\rangle$, 
where $t\mapsto\bigl(\tilde{z}(t),\tilde{\zeta}(t)\bigr)
\in {\cal D}_{\,2\underline{r}}$
is an integral curve for $X_{\tilde{H}}$
and we write $\tilde{\I}(t)$ and $\tilde{\mathbf{\Lambda}}(t)$ 
in place of $I(\tilde{z}(t))$ and $\Lambda(\tilde{z}(t))$, respectively.  
Calculating the derivative, 
using (b) and the fact that $\Lambda$ depends only on
the transverse variable $\tilde \zeta$, we find:
\[ \frac{d\tilde{h}}{dt}=\langle D{h}, X_{\tilde{H}}\rangle=\{{h}, 
\tilde{H}\}
   =\{{h}, \hat{f}_\kappa\}=\langle D{h}, X_{{\hat f}_\kappa}\rangle. \]
Since ${h}(\tilde{z})=\langle\omega^0, I(\tilde{z})\rangle$ depends only 
on $\tilde{z}$ this can be estimated
using only a bound for $\Pi_1 X_{{\hat f}_\kappa}$ which can be obtained in the
same way as (\ref{Cauprep})-(\ref{bmoc}):
$$ {|\Pi_1\,X_{{\hat{f}}_\kappa}|}_{\,5\underline{r}/3}
   \leq\frac{10\,r_2 {|{\hat{f}}_\kappa|}_{\, 2\underline{r}}}{r_1}
   \leq\frac{10\,r_2 2^{-m}\eps}{r_1}
$$
for as long as the solution remains in ${\cal D}_{\,5\underline{r}/3}$,
during which time:
\beq\la{crest}\bigg|\frac{d\tilde{h}}{dt}\bigg|
    =|\langle D{h}, X_{{\hat f}_\kappa}\rangle|\leq\frac{5}{3}\,|\omega^0|\,r_2
   \,{|X_{{\hat{f}}_\kappa}|}_{\,5\underline{r}/3}
   \leq\frac{50\,|\omega^0|r_2^2\,2^{-m}\eps}{3r_1}.
\eeq
Energy conservation,
the convexity assumption (strict positivity of the matrix $A$) and
the coercivity assumption (\ref{spos}) on $\Lambda$ then imply:
\begin{eqnarray}\label{wbent}
   \frac{1}{2M}\,{|\tilde{\I}(t)-{I^0}|}^2
   +\frac{\kappa}{2}\,|\tilde{\zeta}(t)|^2
   & \leq & \frac{1}{2}\,\|A\|\,{|\tilde{\I}(0)-{I^0}|}^2
   +|\tilde{h}(t)-\tilde{h}(0)|+2{|{\hat{g}}|}_{\,5\underline{r}/3}
   +2{|\hat{f}_\kappa|}_{\,5\underline{r}/3}+\kappa\tilde{\mathbf{\Lambda}}(0)
   \nonumber
   \\ & \leq & \frac{1}{2}\,\|A\|\,{|\tilde{\I}(0)-{I^0}|}^2
   +\frac{50\,|\omega^0|r_2^2\,2^{-m}\eps}{3r_1}\,|t|+6\eps
   +\kappa\tilde{\mathbf{\Lambda}}(0). 
\end{eqnarray}
To go further we must relate the initial data in the original and tilde variables.
By (\ref{ass4-2nd}) we know $|\I(0)-{I^0}|\le l_1 r_1<r_1/2$, 
and since we are considering real-valued solutions 
of the Hamiltonian equations, 
\begin{eqnarray*}  
   |z(0)|^2 & = & 2|\I(0)|\le 2(|\I(0)-{I^0}|+|{I^0}|)
   \\ & \le & 2 l_1 r_1+\frac{1}{8}\,r_2^2
   <\Big(\frac{l_1}{2}+\frac{1}{8}\Big)r_2^2\le\frac{r_2^2}{4}\,,
\end{eqnarray*} 
and also $|\zeta(0)|^2\leq 2\mathbf{\Lambda}(0)<(\frac{r_3}{2})^2$ by 
the final conditions in (\ref{ass1-2nd}) and (\ref{ass4-2nd})
of the lemma and (\ref{spos}). Therefore we have 
\[ (z(0),\zeta(0))\in {\cal D}_{\,\underline{r}/2}. \]  
But (a) and (\ref{ass1-2nd})-(\ref{ass2-2nd}) then imply that
\[ |\tilde{z}(0)|\le |z(0)|+\frac{18 m r_2}{r_1}\,\eps T
   \le |z(0)|+\frac{18 l_2r_1}{r_2}
   \le |z(0)|+\frac{9}{2}\,l_2 r_2, \] 
so that for $l_2$ as in (\ref{ass3-2nd}) we get $|\tilde{z}(0)|\le \frac{5r_2}{3}$. 
But then, using (a) again,
\begin{eqnarray}\label{diff0-2nd}
   |\tilde{\I}(0)-{I^0}| & \le & |\tilde{\I}(0)-\I(0)|+|\I(0)-{I^0}|
   \le\frac{1}{2}\,\Big(|\tilde{z}(0)|+|z(0)|\Big)\,|\tilde{z}(0)-\Psi(\tilde{z}(0))|+l_1 r_1
   \nonumber
   \\ & \le & \frac{1}{2}\,\Big(\frac{5r_2}{3}+r_2\Big)\,\frac{18\,mr_2}{r_1}
   \,\eps T+l_1 r_1
   \le\frac{24\,m^2 r_2^2}{r_1}\,\eps T+l_1 r_1\le (24\,l_2+l_1)\,r_1
\end{eqnarray}
by (\ref{ass2-2nd}) and (\ref{ass4-2nd}). 
Thus restricting $|t|$ as in (\ref{conc-2nd}) we obtain from (\ref{wbent}): 
\begin{equation}\label{transI-2nd}
   {|\tilde{\I}(t)-{I^0}|}^2\le M\|A\|\,(24\,l_2+l_1)^2 r_1^2
   +2M(11\,\eps +\kappa\tilde{\mathbf{\Lambda}}(0)). 
\end{equation}
It remains to consider $\tilde{\mathbf{\Lambda}}(0)$. By the fundamental theorem
of calculus and the assumption \eqref{spos2} on $|D\Lambda|$ we have
\[ |\Lambda(\zeta)-\Lambda(\tilde{\zeta})|\leq C_\Lambda(|\zeta|\mu+\mu^2), \] 
since $|\tilde{\zeta}-\zeta|\leq\mu$ by (a). From (\ref{ass1-2nd}) and (\ref{ass2-2nd}) 
it follows that $\mu<r_3$ and thus 
\begin{equation}\label{inter}
   \kappa|\mathbf{\Lambda}(0)-\tilde{\mathbf{\Lambda}}(0)|\leq 2\kappa C_\Lambda r_3\mu
   =\frac{36\kappa C_\Lambda mr_2r_3\eps T}{r_1}
   \leq\frac{36\kappa C_\Lambda mr_1r_2r_3l_0  T}{M}
   \leq\frac{36 l_0  r_1^2}{54 M}
   \leq\frac{r_1^2}{200 M}
\end{equation}
due to $l_0\leq 3/400$.
Hence, using also the final condition in (\ref{ass4-2nd}),
$\kappa\tilde{\mathbf{\Lambda}}(0)\le\frac{r_1^2}{100 M}$
and so by (\ref{transI-2nd}), since the conditions in (\ref{ass3-2nd}) ensure that
$(24\,l_2+l_1)\sqrt{M\|A\|}\le\frac{1}{10}$,
\begin{equation}\label{transIII-2nd}
   {|\tilde{\I}(t)-{I^0}|}^2
   \le\Big(\frac{1}{100}+22\,l_0+\frac{1}{50}\Big) r_1^2<\Big(\frac{r_1}{4}\Big)^2.
\end{equation}
Using the final condition in (\ref{ass1-2nd}) we have similarly from (\ref{wbent}): 
\begin{equation}\label{transII-2nd}
   {|\tilde{\zeta}(t)|}^2
   \le 2\tilde{\mathbf{\Lambda}}(t)\leq
   \frac{\|A\|}{\kappa}\,(24\,l_2+l_1)^2 r_1^2+\frac{22\,\eps}{\kappa}
   +2\tilde{\mathbf{\Lambda}}(0)<\frac{r_1^2}{25\kappa M}<\Big(\frac{r_3}{2}\Big)^2,
\end{equation}
for as long as $\bigl(\tilde{z}(t),\tilde{\zeta}(t)\bigr)
\in {\cal D}_{\,5\underline{r}/3}$
and with $|t|$ restricted as in (\ref{conc-2nd}).
\medskip

Now to deduce (\ref{conc-2nd}) it is necessary to transfer the information
in (\ref{transIII-2nd})--(\ref{transII-2nd}) back to bounds on the original variables
$z$, $\zeta$, $\Lambda$, and $I-I^0$. So let
$t\mapsto (z(t), \zeta(t))\in\R^{2n}\times \R^{2N}$ be the integral curve 
of the original Hamiltonian vector field $X_H$.
Since $(z(0), \zeta(0))\in {\cal D}_{\,\underline{r}/2}$, we can define $t_0>0$ to be 
the longest time such that $(z(t), \zeta(t))\in {\cal D}_{\,\underline{r}}$
for all $|t|\le t_0$, since such a $t_0>0$ exists by continuity.
The point of the following lemma \ref{pobiv-2nd} is to show that
a sufficiently large neighbourhood (to be precise ${\cal D}_{\,\underline{r}}$)
of the approximating periodic orbit determined
by ${I^0}$ is covered by the transformation $\Psi$, 
as a consequence of (a) and the various assumptions on the parameters used.
This ensures that stability information derived for integral curves of 
the transformed Hamiltonian $\tilde H$ does indeed imply stability information
for the integral curves of $H$ on a sufficiently large neighbourhood of 
this periodic orbit. In what follows we write
${\cal D}_{\,\underline{r}}^{({\rm real})}
={\cal D}_{\,\underline{r}}\cap (\R^{2n}\times\R^{2N})$,
and similarly we denote 
\[ 
B_\delta^{({\rm real})}(w, \eta)=\{(z,\zeta)\in\R^{2n}\times \R^{2N}: 
|z-w|+|\zeta-\eta|<\delta\} 
\] 
for $\delta>0$ and $(w,\eta)\in\R^{2n}\times \R^{2N}$.

\begin{lemma}\label{pobiv-2nd}
Under the hypotheses of lemma \ref{lszsk}, $\Psi$ satisfies
$\Psi({\cal D}_{\,5\underline{r}/3}^{({\rm real})})
\supset {\cal D}_{\,\underline{r}}^{({\rm real})}$.
\end{lemma}
{\bf Proof of lemma \ref{pobiv-2nd}\,} According to (a) 
we have ${|\Psi-{\rm id}|}_{\,2\underline{r}}\le\mu=\frac{18\,mr_2}{r_1}\,\eps T$. 
Thus from (\ref{Cauprep}) in conjunction with Cauchy's estimate 
and (\ref{ass1-2nd})--(\ref{ass3-2nd}) we obtain: 
\begin{eqnarray*} 
   {|D\Psi-1|}_{\,5\underline{r}/3}
   & \le &\max\Big\{\frac{3}{r_3},\frac{10\,r_2}{r_1}\,\Big\}
   {|\Psi-{\rm id}|}_{\,2\underline{r}}\\
   &\le &\frac{180\,m r_2^2}{r_1^2}\,\eps T
   \le\frac{180\,m^2 r_2^2}{r_1^2}\,\eps T<180\,l_2<\frac{1}{2}. 
\end{eqnarray*}
Hence $D\Psi(z)$ is invertible for every $z\in {\cal D}_{\,5\underline{r}/3}$, 
and accordingly 
$\Psi: {\cal D}_{\,5\underline{r}/3}\to \Psi({\cal D}_{\,5\underline{r}/3})
=:{\cal W}$ is a real-analytic 
diffeomorphism such that $\|D\Psi^{-1}(w, \eta)\|\le 2$ for 
$(w, \eta)\in {\cal W}$.
Now fix $(w, \eta)\in {\cal D}_{\,\underline{r}}^{({\rm real})}$. 
Then $B_\delta^{({\rm real})}(w, \eta)\subset {\cal D}_{\,3\underline{r}/2}^{({\rm real})}$ 
for $\delta=\frac{r_1}{4\,r_2}<\frac{r_3}{2}$, as can be shown 
using  the first condition in (\ref{ass1-2nd}) and (\ref{Iwz}),
analogously to (\ref{Cauprep}).
Furthermore, for $(w,\eta)\in {\cal D}_{\,\underline{r}}^{({\rm real})}$, 
\[ |(w, \eta)-\Psi(w, \eta)|\le\mu<\frac{\delta}{2} \] 
due to $\frac{18\,mr_2}{r_1}\,\eps T\le\frac{4\cdot 18\,m^2 r_2^2}{r_1^2}\,
\eps T\times\frac{r_1}{4r_2}\le 72\,l_2\delta<\frac{\delta}{2}$. In other words, we have 
$(w, \eta)\in B_{\delta/2}^{({\rm real})}\bigl(\Psi(w,\eta)\bigr)$. 
Next we apply lemma \ref{LipOM} below and use the fact that $\Psi$ is 
real on real vectors,
to deduce that  
\[ \Psi\Big(\,\overline{B_\delta^{({\rm real})}(w, \eta)}\,\Big)
   \supset\overline{B_{\delta/2}^{({\rm real})}\Bigl(\Psi(w, \eta)\Bigr)}. \] 
To summarize, for fixed $(w,\eta)
\in {\cal D}_{\,\underline{r}}^{({\rm real})}$ we obtain 
\[ (w, \eta)\in B_{\delta/2}^{({\rm real})}\Bigl(\Psi(w, \eta)\Bigr)
   \subset\Psi\Big(\,\overline{B_\delta^{({\rm real})}(w, \eta)}\,\Big)
   \subset\Psi\Big(\,\overline{{\cal D}_{\,3\underline{r}/2}^{({\rm real})}}\,\Big)
   \subset\Psi({\cal D}_{\,5\underline{r}/3}^{({\rm real})}), \]  
and this concludes the proof of lemma \ref{pobiv-2nd}. 
\qed\bigskip

\noindent
{\bf Continuation of the proof of lemma \ref{lszsk}}
Due to Lemma \ref{pobiv-2nd}, and referring to the definition of $t_0$, 
we may write $(z(t), \zeta(t))=\Psi(\tilde{z}(t),\tilde{\zeta}(t))$ 
for $|t|\le t_0$ with an integral curve 
$$ t\mapsto\bigl(\tilde{z}(t),\tilde{\zeta}(t)\bigr)
\in {\cal D}_{\,5\underline{r}/3}^{({\rm real})}$$
of $X_{\tilde{H}}$.
Then we can apply (\ref{transIII-2nd})-(\ref{transII-2nd}) to obtain 
$|\tilde{\I}(t)-{I^0}|<r_1/4$ and $|\tilde{\zeta}(t)|<r_3/2$
for $|t|\le\min\{t_\ast, t_0\}$, where $t_\ast=\frac{3\cdot 2^mr_1}{10|\omega^0|r_2^2}$. 
In the same manner as for (\ref{diff0-2nd}) this in turn leads to
\begin{eqnarray}\label{Idiffs-2nd}
   |\I(t)-I^0| & \le & |\I(t)-\tilde{\I}(t)|+|\tilde{\I}(t)-{I^0}|
   \nonumber
   \\ & \le & \frac{1}{2}\,\Big(|\tilde{z}(t)|+|z(t)|\Big)
   \,|\tilde{z}(t)-\Psi(\tilde{z}(t))|+\frac{r_1}{4}
   \nonumber
   \\ & \le & \Big(24\,l_2+\frac{1}{4}\Big) r_1<\frac{r_1}{2}
\end{eqnarray}
for $|t|\le\min\{t_\ast, t_0\}$ and as
$24l_2<\frac{1}{4}$. Since also $r_1<r_2^2/4$, 
this implies that for such times
\begin{equation}\label{d1s}
   |z(t)|^2 =  2|\I(t)|\le 2\Big(|\I(t)-{I^0}|+|{I^0}|\Big)
   \le2\Big(\frac{r_2^2}{8}
   +\frac{r_2^2}{16}\Big)<r_2^2.
\end{equation}
Also, as in the derivation of (\ref{inter}) and by (\ref{transII-2nd}), 
we have for $|t|\le\min\{t_\ast, t_0\}$
\begin{equation}\label{d3s}
   \kappa\mathbf{\Lambda}(t)\leq\kappa |\mathbf{\Lambda}(t)-\tilde{\mathbf{\Lambda}}(t)|
   +\kappa\tilde{\mathbf{\Lambda}}(t)
   \leq 2C_\Lambda r_3\kappa\mu+\kappa\tilde{\mathbf{\Lambda}}(t)
   \leq\frac{r_1^2}{200 M}+\frac{r_1^2}{50 M}
   <\frac{r_1^2}{16 M},
\end{equation}
and furthermore by (a),
\begin{equation}\label{d2s}
|\zeta(t)|\leq |\tilde{\zeta}(t)-\zeta(t)|+|\tilde{\zeta}(t)|
\leq\mu+\frac{r_3}{2}<r_3, 
\end{equation}
the latter since $\mu\le 18 l_2 r_1/r_2<36 l_2 r_3<r_3/2$.
Altogether from (\ref{d1s}) and (\ref{d2s}) 
we conclude that $\min\{t_\ast, t_0\}<t_0$, 
or in other words $\min\{t_\ast, t_0\}=t_\ast$,
and so the assertions in (\ref{conc-2nd}) follow as a consequence 
of (\ref{Idiffs-2nd}), (\ref{d3s}), and (\ref{d2s}). \qed\bigskip

\noindent
{\bf Completion of proof of theorem \ref{gszsk}\,}
We now aim to show that the stability bound (\ref{conc-2nd}), 
applied in the neighbourhood of the approximating periodic orbit obtained
prior to lemma \ref{lszsk}, implies (\ref{gszslk}).
Recall that the period of the approximating periodic orbit is $T=\tau/\theta^2$,
and define $r_2=8\theta$ and
\begin{equation}\label{impdef-2nd}
   m=\delta\,[\theta^{-a}],\quad
   r_1=\frac{L\theta^{2+a}}{\tau},\quad 
   \mbox{and}\quad r_3=P\theta^{1+2a(1-n)},
\end{equation}
where $\delta$, $L$, $P>0$ will be fixed below.
To ensure that ${|f_\kappa|}_{\,3\underline{r}}\le\eps$
define $\eps=C_1\theta^5$, so that
\[ {|f_\kappa|}_{\,3\underline{r}}\le\sup\,\{|f_\kappa(z, \zeta)|: |z|\le 3r_2, |\zeta|\le r_3\}
   \le C_0\bigl[(3r_2)^5+r_3^2(3r_2)^4+\kappa r_3^2(3r_2)\bigr]
   \le\eps\,, \] 
with the choice $C_1=(24^5+24^4 P^2+24P^2)C_0$; 
here we have used $\theta<1$, $\kappa\le\theta^{2+4a(n-1)}$ due to (\ref{agsz}), 
and the restriction $a<\frac{1}{4(n-1)}$ from the beginning of the theorem statement.
We will now verify that having defined $l_0, l_1, l_2>0$ by (\ref{ass3-2nd}), 
the conditions (\ref{ass1-2nd})--(\ref{ass4-2nd}) can be made to hold
by making $\theta$ sufficiently small and choosing $\delta$, $L$, $P$ 
appropriately. 

{\em The conditions in (\ref{ass1-2nd})}. 
To start with
\[ \frac{r_1}{r_2^2}=\frac{L\theta^a}{64\tau}\le\frac{L\theta^a}{64\pi} 
\quad\hbox{and}\quad
\frac{r_1}{r_2 r_3}\le\frac{L\theta^{a(2n-1)}}{8\pi P} 
\] 
by (ii), and hence the first condition of (\ref{ass1-2nd}) 
holds if $\theta$ is small enough (depending upon $L$, $P$, $a$, $n$). Next, 
\[ \frac{\eps}{r_1^2}=\frac{C_1\theta^5\tau^2}{L^2\theta^{4+2a}}
   \le\frac{16\pi^2 C_1}{L^2}\,\theta^{1-2an} \] 
by (ii) shows that we can fulfil the second condition of (\ref{ass1-2nd}) 
for $\theta$ sufficiently small (depending upon $L$, $C_0$, $P$, $a$, $n$), 
due to $a<\frac{1}{2n}$. 
Concerning the condition on $|{I^0}|$ in (\ref{ass1-2nd}), here 
\[ |{I^0}|\le |\I(0)-{I^0}|+|\I(0)|\le nC_A\,\frac{\theta^{2+a}}{\tau}+\theta^2
   \le\Big(nC_A\,\frac{\theta^a}{\pi}+1\Big)\theta^2\le 2\theta^2 \]
by (i) and (ii) for $\theta$ small enough (depending upon $C_A$, $a$, $n$). Hence 
\[ \frac{|{I^0}|}{r_2^2}\le \frac{2\theta^2}{64\theta^2}<\frac{1}{16}. \]
The final condition in (\ref{ass1-2nd}) reads as
\[ 1\ge\frac{r_1^2}{4\kappa Mr_3^2}=\frac{L^2}{4MP^2\tau^2} \] 
recall (\ref{agsz}). Since $\tau\ge\pi$, this follows from 
\begin{equation}\label{delL4-2nd} 
   \frac{L^2}{P^2}\le 4\pi^2 M.
\end{equation} 
 
{\em The conditions in (\ref{ass2-2nd})}. 
Next
\[ m r_1 T=\delta\,[\theta^{-a}]\,\frac{L\theta^{2+a}}{\tau}\,\frac{\tau}{\theta^2}\le\delta L \] 
whence the restriction
\begin{equation}\label{delL1-2nd} 
   \delta L\le\frac{1}{324\,\|A\|}
\end{equation} 
is sufficient to validate the first condition of (\ref{ass2-2nd}).
For the second condition of (\ref{ass2-2nd}) calculate
\[ \frac{m^2\eps T r_2^2}{r_1^2}
   \le\frac{\delta^2\theta^{-2a}C_1\theta^5\tau^3\times 64\theta^2}{\theta^2 L^2\theta^{4+2a}}
   =\frac{64\,C_1\delta^2}{L^2}\,\theta^{1-4a}\,\tau^3
   \le\frac{4096\,\pi^3C_1\delta^2}{L^2}\,\theta^{1-a(1+3n)} \] 
by (ii). Since $a<\frac{1}{1+3n}$, the right-hand side is smaller 
than $l_2$ from (\ref{ass3-2nd}), if $\theta$ is sufficiently small
(depending upon $\delta$, $L$, $C_0$, $P$, $a$, $n$). 
The final condition in (\ref{ass2-2nd}) is 
\[ 1\le\frac{r_1}{54 C_\Lambda\kappa m T r_2 r_3}
   =\frac{L\theta^{a(1-2n)}}{432 C_\Lambda\delta P\tau^2\,[\theta^{-a}]}, \] 
which due to (ii) is a consequence of $432 C_\Lambda\delta P (4\pi)^2 [\theta^{-a}]\le L\theta^{-a}$. 
This in turn holds if
\begin{equation}\label{delL3-2nd}
   \delta P\leq\frac{L}{6912\pi^2 C_\Lambda}.
\end{equation}

{\em The conditions in (\ref{ass4-2nd})}. The first one holds because 
\[ \frac{|\I(0)-{I^0}|}{r_1}\le nC_A\,\frac{\theta^{2+a}}{\tau}\,\frac{\tau}{L\theta^{2+a}}
   =\frac{nC_A}{L} \] 
by (i). Accordingly, we need to have 
\begin{equation}\label{delL2-2nd} 
   \frac{nC_A}{L}\le l_1
\end{equation} 
for $l_1$ from (\ref{ass3-2nd}). The second condition holds because
$r_1=\frac{L\theta^{2+a}}{\tau}\geq\frac{L\theta^{2+an}}{4\pi}$ by (ii), 
so that due to (\ref{agsz}) 
\[ \kappa\mathbf{\Lambda}(0)\le C_E\theta^{4+2an}\leq \frac{C_E(4\pi)^2r_1^2}{L^2}=\frac{r_1^2}{200M},
   \quad\hbox{taking}\quad C_E=\frac{L^2}{(4\pi)^2200 M}. \] 

Altogether, the conditions necessary to apply lemma \ref{lszsk} will be satisfied provided that
the restrictions in (\ref{delL4-2nd}), (\ref{delL1-2nd}), (\ref{delL3-2nd}) and (\ref{delL2-2nd}) hold.
The latter can be achieved by explicitly taking 
\[ L=\frac{nC_A}{l_1},\quad P=\frac{L}{2\pi\sqrt{M}}
   \quad\mbox{and}\quad
   \delta=\min\Big\{\frac{1}{324\,\|A\| L}\,,\,\frac{L}{6912\pi^2 C_\Lambda P}\Big\}. \] 
Therefore lemma \ref{lszsk} can be used, and we deduce from (\ref{conc-2nd}) that 
\[ |\I(t)-{I^0}|\le r_1,\quad\kappa\mathbf{\Lambda}(t)\le\frac{r_1^2}{16M}, 
   \quad\mbox{and}\quad |\zeta(t)|\leq r_3\quad\hbox{for}\quad 
   |t|\le\frac{3\times 2^m r_1}{10|\omega^0|r_2^2}=:t_\ast. \] 
Now combine the former with (i) to bound
$$
|\I(t)-\I(0)|\leq |\I(t)-{I^0}|+|{I^0}-\I(0)|\leq r_1+nC_A
\frac{\theta^{2+a}}{\tau}\le\frac{2L}{\pi}\,\theta^{2+a}
$$
(since clearly $L>nC_A$) for $|t|\le t_\ast$. Thus, recalling that
$\tau\ge\pi$  by (ii), we can achieve the bounds in \eqref{gszslk} with 
\[ K=\max\Big\{\frac{2L}{\pi}, \frac{L^2}{16 M\pi^2}\Big\}\,, 
\] 
and
\[ t_\ast=\frac{3\times 2^m L\theta^a}{640\,|\omega^0|\tau}
   \ge\frac{3\times 2^m L\theta^{an}}{2560\,\pi |\omega^0|}\,. 
\] 
It follows that with $B=\frac{3L}{2560\,\pi |\omega^0|}$ 
and $m$ as in (\ref{impdef-2nd}),
\[ \ln t_\ast = \ln B+an\ln\theta+(\ln 2)\delta\,[\theta^{-a}]
   \ge k\,\theta^{-a}\,,\] 
for $\theta$ small enough and any
$k<(\ln 2)\delta$, thus completing the proof of (\ref{gszslk}).
\qed\bigskip

%%%%%%%%%%%%%%%%%%%%%%%%%%%%%%%%%%%%%%%%%%%%%%%%%%%%%%%%%%%%%%%%%%%%%%%%%%%%%%%%%%%%%%%%

\setcounter{equation}{0}

\section{A variant of the main theorem}

In some applications it may be desirable to prove that the stability estimates
hold for sufficiently small $\kappa$ in an open set in $\R^{2n}$ which is 
essentially determined by the unperturbed $z$ motion. In these circumstances
the following variant of theorem \ref{gszsk} is natural:

\begin{theorem}\label{gszskvar} 
Let $t\mapsto\bigl(z(t),\zeta(t)\bigr)$  be
an integral curve of the real-analytic Hamiltonian $H$ verifying
\eqref{ra}, \eqref{spos}-\eqref{spos2} and
\begin{equation}\label{tenth}
   |f_\kappa(z, \zeta)|\leq C_0\bigl(|z|^5+\kappa |\zeta|^2|z|\bigr)
\end{equation}
(in place of \eqref{fifth}).
Fix $a\in ]0, \min\{\frac{1}{4(n-1)}, \frac{1}{1+3n}\}[$. 
Then there exist positive numbers $C_E$, $\theta_0<1$, $K$, $k$ 
with the following properties. If the initial data are such that 
$\I(0)=I(z(0))$ and $\mathbf{\Lambda}(0)=\Lambda(\zeta(0))$ satisfy
\begin{equation}\label{agszv}
   |\I(0)|\leq\theta^2\,,\quad
   \kappa\mathbf{\Lambda}(0)\leq C_E\theta^{4+2an}
   \quad\hbox{and}\quad
   0<\kappa \leq \theta^{2+2a(2n-1)}\,,
\end{equation}
for $0<\theta\le\theta_0$, 
then $\I(t)=I(z(t))$ and $\mathbf{\Lambda}(t)=\Lambda(z(t))$ satisfy
the stability estimate (\ref{gszslk}) as in theorem \ref{gszsk}.
\end{theorem}
\noindent
{\bf Proof of theorem \ref{gszskvar}\,}
Only a small modification of the 
proof of theorem \ref{gszsk} is needed. We start by using periodic
approximation and lemma \ref{lszsk} in identical fashion, but then
in \eqref{impdef-2nd} replace the definition of $r_3$ by
\beq\la{impdef-3rd}
r_3=\frac{P\theta^{2+a}}{\sqrt{\kappa}}
\eeq
(leaving the definitions of $r_1,r_2,m$ unchanged). When $\kappa$ is
equal to its maximum allowed value, $\theta^{2+2a(2n-1)}$, this reproduces 
the value of $r_3$
in \eqref{impdef-2nd}, but as $\kappa$ gets smaller $r_3$ defined in 
\eqref{impdef-3rd} increases in such a way that $\kappa r_3^2$ is
unchanged. With this understood it is easy to see that the
conditions \eqref{ass2-2nd}-\eqref{ass4-2nd} in lemma \ref{lszsk}
continue to hold with the new definition of $r_3$, and hence 
the estimate (\ref{gszslk}) holds as a consequence of that lemma
exactly as in the completion of the proof of theorem \ref{gszsk}.
\qed\bigskip

%%%%%%%%%%%%%%%%%%%%%%%%%%%%%%%%%%%%%%%%%%%%%%%%%%%%%%%%%%%%%%%%%%%%%%%%%%%%%%

\setcounter{equation}{0}

\section{Some auxiliary results}

In this section ${|\cdot|}_\infty$ denotes the maximum norm on $\R^n$,
i.e.~${|x|}_\infty=\max_{1\leq j\leq n}|x_j|$.

\begin{lemma}[Dirichlet]\label{diri-lem} For every $Q\in\N$ and $\omega\in\R^n$,
\[ \min_{q\in\{1, \ldots, Q\}}\,\min_{p\in\Z^n}\,{|q\omega-p|}_\infty\le\frac{1}{Q^{1/n}}. \]
\end{lemma}
{\bf Proof\,} See \cite[Thm.~1B, p.~34]{schmidt}. \qed\bigskip

\begin{cor}\label{dircor}
For every $Q\in\N$ and $\omega\in\R^n$
such that ${|\omega|}_\infty>1$
there exists $\omega^0\in\R^n$ and 
$T\in [2\pi (1-{|\omega|}_\infty^{-1})\,,\,2\pi Q]$ 
and 
such that $\omega^0 T\in 2\pi\Z^n$ and
\begin{equation}\label{diffom}
   {|\omega-\omega^0|}_\infty\le\frac{2\pi}{T Q^{1/(n-1)}}.
\end{equation}
\end{cor}
{\bf Proof\,} We may assume that $\omega_n={|\omega|}_\infty>1$,
where $\omega=(\omega_1, \ldots, \omega_n)$. Then write
\[ \omega=\frac{\omega_n}{[\omega_n]}\Big(\frac{[\omega_n]}{\omega_n}\,\omega_1,
   \ldots, \frac{[\omega_n]}{\omega_n}\,\omega_{n-1}, [\omega_n]\Big) \]
and apply Lemma \ref{diri-lem} to find $q\in\{1, \ldots, Q\}$ and $p\in\Z^{n-1}$
such that
\begin{equation}\label{dirom}
   \Big|q\,\frac{[\omega_n]}{\omega_n}\,\omega_j-p_j\Big|\le\frac{1}{Q^{1/(n-1)}},
   \quad 1\le j\le n-1.
\end{equation}
Defining $T=2\pi q\,[\omega_n]/\omega_n$ and
\[ \omega^0=\frac{\omega_n}{[\omega_n]}\Big(\frac{p_1}{q},
   \ldots, \frac{p_{n-1}}{q}, [\omega_n]\Big), \]
it follows that $T\omega^0=2\pi (p_1, \ldots, p_{n-1}, q[\omega_n])\in 2\pi\Z^n$.
Furthermore, $1-1/\omega_n\le [\omega_n]/\omega_n\le 1$ and $1\le q\le Q$ yield
the bound on $T$. Finally, (\ref{diffom}) is a direct consequence of (\ref{dirom}).
\qed\bigskip

As a further corollary we obtain the density of periodic orbits in sufficiently
small neighbourhoods of elliptic equilibria for convex integrable Hamiltonians.
The frequency $\omega^0$ is called $T-${\it periodic} if
$\omega^0 T\in2\pi\Z^n$.
Now to be precise consider a real analytic Hamiltonian on $\R^{2n}$ of the form
\begin{equation}\label{H0def}
\langle\alpha, I\rangle+\frac{1}{2}\,\langle AI, I\rangle+g(I)
\end{equation} 
where $A$ is a strictly positive matrix and $g(I)=O(|I|^3)$, and the notation
is as in the introduction. The function
\begin{equation}\label{fmap}
I\mapsto\Omega(I)
=\alpha+AI+Dg(I)=\alpha+AI+O(|I|^2) 
\end{equation}
is invertible in a neighbourhood of the origin in 
$\R^n$ by the inverse function theorem,
since $D\Omega(I)=A+O(|I|)$ is invertible for $|I|$ small enough.
The smooth inverse $\Omega^{-1}$ is defined on a 
neighbourhood of $\alpha=\Omega(0)$.

\begin{cor}\label{pds} Given a function $\Omega:\R^n\to\R^n$ 
as in (\ref{fmap}), with $\alpha\in\R^n\setminus\{0\}$, and a number $a>0$,
there exist $C>0$ (depending upon $\Omega$) and $\theta_0>0$
(depending upon $\Omega,\alpha,a,n$) such that the following holds:
if $I\in\R^n$ and $|I|=\theta^2$ for some $0<\theta\le\theta_0$,
then there exist ${I^0}\in\R^n$ and $\tau>0$ such that
\begin{enumerate}
\item
${|I-{I^0}|}_\infty\leq C\,\frac{\theta^{2+a}}{{\tau}}$, 
\item
$\pi\le 2\pi (1-\theta^2{|\Omega(I)|}_\infty^{-1})\le {\tau}\le 4\pi\theta^{-a(n-1)}$, and 
\item
$\omega^0=\Omega({I^0})$ is ${\tau}/\theta^2$-periodic, 
i.e.~$\omega^0\frac{{\tau}}{\theta^2}\in 2\pi\Z^n$.
\end{enumerate}
\end{cor}
{\bf Proof\,} By the above remarks there is $\eps>0$ such that
$\Omega: B_\eps(0)\to\Omega(B_\eps(0))=:U$ is smoothly invertible
and $C=2\pi {\|D\Omega^{-1}\|}_{L^\infty(\overline{U})}<\infty$.
Choose $\delta$ obeying $0<\delta<|\alpha|_\infty$ 
such that $B_\delta(\alpha)\subset U$.
Next fix $\theta_0>0$ sufficiently small that for $0<\theta\le\theta_0$ 
and $|I|=\theta^2$ there holds:
\[ {|\Omega(I)-\alpha|}_\infty<\delta/2,
   \quad\theta^2<\min\{\eps, {|\alpha|}_\infty/4\},\quad\theta^{a(n-1)}<1,
   \quad 2\theta^{2+a}<\delta/2; \]
hence $\theta_0$ depends on $\Omega,\alpha,a$ and $n$. Now if $|I|=\theta^2$
for some $0<\theta\le\theta_0$, then $|I|<\eps$ 
ensures that $\omega=\Omega(I)\in U$ is well-defined, and
$|\omega|_\infty>|\alpha|_\infty/2$ since $\delta<|\alpha|_\infty$.
Putting $Q=[\theta^{-a(n-1)}]+1$ and $\tilde{\omega}=\theta^{-2}\omega$,
we have ${|\tilde{\omega}|}_\infty\ge\theta^{-2}{|\alpha|}_\infty/2>2>1$.
Therefore corollary \ref{dircor} applies and yields the existence of
${\tau}>0$ and $\tilde{\omega}^0\in\R^n$ such that $\tilde{\omega}^0 {\tau}\in 2\pi\Z^n$,
$2\pi (1-{|\tilde{\omega}|}_\infty^{-1})\le {\tau}\le 2\pi Q$ and
\begin{equation}\label{tilddiff}
   {|\tilde{\omega}-\tilde{\omega}^0|}_\infty\le\frac{2\pi}{{\tau} Q^{1/(n-1)}}.
\end{equation}
Also, ${\tau}\ge 2\pi (1-{|\tilde{\omega}|}_\infty^{-1})\ge\pi$ follows
since ${|\tilde{\omega}|}_\infty>2$.
Defining $\omega^0=\theta^2\tilde{\omega}^0$, we get $\omega^0\frac{{\tau}}{\theta^2}\in 2\pi\Z^n$.
Furthermore,
\[ {|\omega-\omega^0|}_\infty=\theta^2 {|\tilde{\omega}-\tilde{\omega}^0|}_\infty
   \le\frac{2\pi\theta^2}{{\tau} Q^{1/(n-1)}}\le\frac{2\theta^2}{{\theta}^{-a}}=2\theta^{2+a}
   <\delta/2 \]
implies that ${|\omega^0-\alpha|}_\infty\le {|\omega^0-\omega|}_\infty+{|\Omega(I)-\alpha|}_\infty
<\delta/2+\delta/2=\delta$, and consequently $\omega^0\in U$ so that ${I^0}=\Omega^{-1}(\omega^0)$
is well defined. Then (ii) follows from $2\pi Q\le 2\pi (\theta^{-a(n-1)}+1)\le 4\pi\theta^{-a(n-1)}$.
Finally, concerning (i) it suffices to note that since
both $\omega$ and $\omega^0$ lie in the ball $B_\delta(\alpha)\subset U$
\[ {|I-{I^0}|}_\infty={|\Omega^{-1}(\omega)-\Omega^{-1}(\omega^0)|}_\infty
   \le {\|D\Omega^{-1}\|}_{L^\infty(\overline{U})}\,{|\omega-\omega^0|}_\infty
   \le {\|D\Omega^{-1}\|}_{L^\infty(\overline{U})}\,\frac{2\pi\theta^2}{{\tau}{\theta}^{-a}}
   =C\,\frac{\theta^{2+a}}{{\tau}} \]
by (\ref{tilddiff}), completing the proof.
\qed\bigskip

We also need the following quantitative version of the inverse mapping theorem.

\begin{lemma}\label{LipOM} Let $X, Y$ be Banach spaces and suppose that
$U\subset X$ is open. If $\Psi: U\to\Psi(U)\subset Y$
is a homeomorphism, $\Psi^{-1}$ is Lipschitz continuous
with constant 
${\rm Lip}(\Psi^{-1})<\lambda$, and $\overline{B_r(x)}\subset U$, then
\[ \Psi(\overline{B_r(x)})\supset\overline{B_{r/\lambda}(\Psi(x))}. \]
\end{lemma}
{\bf Proof\,} See \cite[Prop.~I.3, p.~50]{shub}.
\qed\bigskip

%%%%%%%%%%%%%%%%%%%%%%%%%%%%%%%%%%%%%%%%%%%%%%%%%%%%%%%%%%%%%%%%%%%%%%%%%%%%%%%%%%%%%%%%%%%%%

\end{document}